\documentclass[11pt]{article}%
\usepackage{amsmath}
\usepackage{amssymb}
\usepackage{amsfonts}
\usepackage{graphicx}%
\makeatletter
\newfont{\footsc}{cmcsc10 at 8truept}
\newfont{\footbf}{cmbx10 at 8truept}
\newfont{\footrm}{cmr10 at 10truept}
\newtheorem{theorem}{Theorem}[section]

\newtheorem{lemma}[theorem]{Lemma}

\begin{document}

\title{Some new results in extremal graph theory}
\author{Vladimir Nikiforov \thanks{Research supported by NSF Grant DMS-0906634.}\\{\small Department of Mathematical Sciences, University of Memphis, Memphis,
TN 38152}\\{\small e-mail:} {\small vnikifrv@memphis.edu}}
\maketitle

\begin{abstract}
In recent years several classical results in extremal graph theory have been
improved in a uniform way and their proofs have been simplified and
streamlined. These results include a new Erd\H{o}s-Stone-Bollob\'{a}s theorem,
several stability theorems, several saturation results and bounds for the
number of graphs with large forbidden subgraphs.

Another recent trend is the expansion of spectral extremal graph theory, in
which extremal properties of graphs are studied by means of eigenvalues of
various matrices. One particular achievement in this area is the casting of
the central results above in spectral terms, often with additional
enhancement. In addition, new, specific spectral results were found that have
no conventional analogs.

All of the above material is scattered throughout various journals, and since
it may be of some interest, the purpose of this survey is to present the best
of these results in a uniform, structured setting, together with some
discussions of the underpinning ideas. \medskip

\newpage

\end{abstract}
\tableofcontents

\newpage

\section{Introduction}

The purpose of this survey is to give a systematic account of two recent lines
of research in extremal graph theory. The first one, developed in
\cite{BoNi04},\cite{BoNi09b},\cite{BoNi10},\cite{Nik08,Nik10c}, improves a
number of classical results grouped around the theorem of Tur\'{a}n. The main
progress is along the following three guidelines: replacing fixed parameters
by variable ones; giving explicit conditions for the validity of the
statements; developing and using tools of general scope. Among the results
obtained are a new Erd\H{o}s-Stone-Bollob\'{a}s theorem (see Section
\ref{ErSt}), several stability theorems (see Section \ref{stab}), several
saturation results, and bounds for the number of graphs without given large subgraphs.

The second line of research, developed in \cite{BoNi07},\cite{Nik02,Nik10b},
can be called \emph{spectral extremal graph theory}, where connections are
sought between graph properties and the eigenvalues of certain matrices
associated with graphs. As a result of this research, much of classical
extremal graph theory has been translated into spectral statements, and this
translation has also brought enhancement. Among the results obtained are
spectral forms of the Tur\'{a}n theorem and the Erd\H{o}s-Stone-Bollob\'{a}s
theorem, several stability theorems, along with new bounds for the
Zarankiewicz problem (What is the maximum number of edges in a graph with no
$K_{s,t}$?).

In the course of this work a few tools were developed, which help to cast
systematically some classical results and their proofs into spectral form. The
use of this machinery is best exhibited in \cite{Nik09c}, where we gave a new
stability theorem and also its spectral analog - Theorems \ref{tStab} and
\ref{stStab} below. As an illustration, in Section \ref{spf} we outline the
proofs of these two results.

We believe that ultimately the spectral approach to extremal graph theory will
turn out to be more fruitful than the conventional one, albeit it is also more
difficult, and is still underdeveloped. Indeed, most statements in
conventional terms can be cast and proved in spectral terms, but in addition
to that, there are a lot of specific spectral results (say, Theorem
\ref{minK}) with no conceivable conventional setting.

The rest of the survey is organized as follows. To keep the beginning
straightforward, the bulk of the necessary notation and the basic facts have
been shifted to Section \ref{not}, although some definitions are given also
where appropriate. Section \ref{class} covers the conventional, nonspectral
problems, while Section \ref{spec} presents the spectral results. In Section
\ref{tool}, we have collected some basic and more widely applicable
statements, which we have found useful on more than one occasion. Finally
Section \ref{spf} presents some proof techniques for illustration, and in fact
these are the only proofs in this survey.

\section{\label{class}New results on classical extremal graph problems}

In extremal graph theory one investigates how graph properties depend on the
value of various graph parameters. In a sense almost all of graph theory deals
with extremal problems, but there is a bundle of results grouped around
Tur\'{a}n's theorem \cite{Tur41}, that undoubtedly constitutes the core of
extremal graph theory. To state this celebrated theorem, which has stimulated
researchers for more than six decades, recall that for $n\geq r\geq2,$ the
Tur\'{a}n graph $T_{r}\left(  n\right)  $ is the complete $r$-partite graph of
order $n$ whose class sizes differ by at most one. We let $t_{r}\left(
n\right)  =e\left(  T_{r}\left(  n\right)  \right)  .$

\begin{theorem}
If $G$ is a graph of order $n,$ with no complete subgraph of order $r+1,$ then
$e\left(  G\right)  \leq t_{r}\left(  n\right)  $ with equality holding only
when $G=T_{r}\left(  n\right)  .$
\end{theorem}

Here is a more popular, but slightly weaker version, which we shall call the
\emph{concise Tur\'{a}n theorem:}\medskip

\emph{If }$G$\emph{ is a graph of order }$n,$\emph{ with }$e\left(  G\right)
>\left(  1-1/r\right)  n^{2}/2,$ \emph{then }$G$\emph{ contains a complete
subgraph of order }$r+1.$\medskip

No doubt, Tur\'{a}n's theorem is a nice combinatorial statement and it is not
too difficult to prove as well. However, its external simplicity is
incomparable with its real importance, since this theorem is a cornerstone on
which rest much more general statements about graphs. Thus, in this survey, we
shall meet the Tur\'{a}n graph $T_{r}\left(  n\right)  $ and the numbers
$t_{r}\left(  n\right)  $ on numerous occasions.

\subsection{The extremal problems that are studied}

Among the many questions motivated by Tur\'{a}n's theorem, the ones that we
will discuss in Section \ref{class} fall into the following three broad
classes:\medskip

\emph{(1) Which subgraphs are present in a graph }$G$\emph{ of order }%
$n$\emph{ whenever }$e\left(  G\right)  >t_{r}\left(  n\right)  \emph{\ and}$
$n$ \emph{is sufficiently large?}\smallskip

As we shall see, here the range of $e\left(  G\right)  -t_{r}\left(  n\right)
$ determines different problems: when $e\left(  G\right)  -t_{r}\left(
n\right)  =o\left(  n\right)  $ we have \emph{saturation problems}, and when
$e\left(  G\right)  -t_{r}\left(  n\right)  =o\left(  n^{2}\right)  $, we have
\emph{Erd\H{o}s-Stone type problems}.\smallskip

Other questions that we will be interested in give rise to the so called
\emph{stability problems, }concerning near-maximal graphs without forbidden
subgraphs.\smallskip

\emph{(2) Suppose that }$H_{n}$\emph{ is a graph which is present in any graph
}$G$\emph{ of order }$n$\emph{ whenever }$e\left(  G\right)  >t_{r}\left(
n\right)  ,$\emph{ but }$H_{n}$\emph{ is not a subgraph of the Tur\'{a}n graph
}$T_{r}\left(  n\right)  $\emph{. We can ask the following questions:}\medskip

-\emph{ What can be the structure of an }$H_{n}$\emph{-free graph }$G$\emph{
of order }$n$ \emph{if }$e\left(  G\right)  >t_{r}\left(  n\right)  -f\left(
n\right)  ,$ \emph{where }$f\left(  n\right)  \geq0$\emph{ and} $f\left(
n\right)  =o\left(  n^{2}\right)  ?$\medskip

-\emph{ What can be the structure of an }$H_{n}$\emph{-free graph }$G$\emph{
of order }$n,$ \emph{with minimum degree }$\delta\left(  G\right)  >\left(
1-c\right)  \delta\left(  t_{r}\left(  n\right)  \right)  $\emph{ for some
sufficiently small }$c>0?$\medskip

Obviously these two general questions have lots of variations, many of which
are intensively studied due to their applicability in other extremal
problems.\medskip

Finally, recall that a long series of results deals with the number of graphs
having some monotone or hereditary properties. Here we will discuss a similar
and natural question which, however, goes beyond this paradigm:\medskip

\emph{(3) Let }$\left\{  H_{n}\right\}  $\emph{ be a sequence of graphs with
}$v\left(  H_{n}\right)  =o\left(  \log n\right)  .$ \emph{How many }$H_{n}%
$\emph{-free graphs of order }$n$ \emph{are there}$?$

\subsection{\label{ErSt}Erd\H{o}s-Stone type problems}

We write $K_{r}\left(  s_{1},...,s_{r}\right)  $ for the complete $r$-partite
graph with class sizes $s_{1},...,s_{r},$ and set for short%
\[
K_{r}\left(  p\right)  =K_{r}\left(  p,...,p\right)  \text{ \ \ and \ \ }%
K_{r}\left(  p;q\right)  =K_{r}\left(  p,...,p,q\right)  .
\]
Let us recall the fundamental theorem of Erd\H{o}s and Stone \cite{ErSt46}.

\begin{theorem}
\label{ES}For all $c>0$ and natural $r,p,$ there is an integer $n_{0}\left(
p,r,c\right)  $ such that if $G$ is a graph of order $n>n_{0}\left(
p,r,c\right)  $ and $e\left(  G\right)  \geq\left(  1-1/r+c\right)  n^{2}/2,$
then $G$ contains a $K_{r+1}\left(  p\right)  $.
\end{theorem}

Noting that $t_{r}\left(  n\right)  \thickapprox\left(  1-1/r\right)
n^{2}/2,$ we see the close relation of Theorem \ref{ES} to Tur\'{a}n's
theorem. In fact, Theorem \ref{ES} answers a fairly general question:
\emph{what is the maximum number of edges }$e\left(  n,H\right)  $ \emph{in a
graph of order }$n$ \emph{that does not contain a fixed }$\left(  r+1\right)
$\emph{-chromatic subgraph }$H?$ Theorem \ref{ES} immediately implies that
$e\left(  n,H\right)  \leq\left(  1-1/r+o\left(  1\right)  \right)  n^{2}/2.$
On the other hand, $T_{r}\left(  n\right)  $ contains no $\left(  r+1\right)
$-chromatic subgraphs, and so, $e\left(  n,H\right)  =\left(  1-1/r+o\left(
1\right)  \right)  n^{2}/2.$

Write $g\left(  n,r,c\right)  $ for the maximal $p$ such that every graph $G$
of order $n$ with%
\[
e\left(  G\right)  \geq\left(  1-1/r+c\right)  n^{2}/2
\]
contains a $K_{r+1}\left(  p\right)  .$ For almost 30 years the order of
magnitude of $g\left(  n,r,c\right)  $ remained unknown; it was established
first by Bollob\'{a}s and Erd\H{o}s in \cite{BoEr73}, as given below. This
simplest quantitative form of the Erd\H{o}s-Stone theorem we call the
\emph{Erd\H{o}s-Stone-Bollob\'{a}s theorem}.

\begin{theorem}
\label{BE}There are constants $c_{1},c_{2}>0$ such that
\[
c_{1}\log n\leq g\left(  n,r,c\right)  \leq c_{2}\log n.
\]

\end{theorem}

Subsequently the function $g\left(  n,r,c\right)  $ was determined with great
precision in \cite{BES76}, \cite{ChSz83}, \cite{BoKo94}, \cite{Ish02}, to name
a few milestones. However, since Szemer\'{e}di's Regularity Lemma is a
standard tool in this research, the results are confined to fixed $c,$ and $n$
extremely large.

To overcome these restrictions, in \cite{Nik08}, we proposed a different
approach, based on the expectation that the presence of many copies of a given
subgraph $H$ must imply the existence of large blow-ups of $H.$ As a
by-product, this approach gave results in other directions as well, which
otherwise do not seem too close to the Erd\H{o}s-Stone theory; two such topics
are outlined in \ref{copy} and \ref{r-gr}.

\subsubsection{Refining the Erd\H{o}s-Stone-Bollob\'{a}s theorem}

The general idea above is substantiated for cliques in the following two
theorems, given in \cite{Nik08}.

\begin{theorem}
\label{Nik1} Let $r\geq2,$ let $c$ and $n$ be such that
\[
0<c<1/r!\text{ \ \ and \ \ }n\geq\exp\left(  c^{-r}\right)  ,
\]
and let $G$ be a graph of order $n.$ If $k_{r}\left(  G\right)  >cn^{r},$ then
$G$ contains a $K_{r}\left(  s;t\right)  $ with $s=\left\lfloor c^{r}\log
n\right\rfloor $ and $t>n^{1-c^{r-1}}.$
\end{theorem}

In a nutshell, Theorem \ref{Nik1} says that if a graph contains many
$r$-cliques, then it has large complete $r$-partite subgraphs. Hence, to
obtain Theorem \ref{BE}, all we need to prove is that the hypothesis of the
Erd\H{o}s-Stone theorem implies the existence of sufficiently many
$r$-cliques. This implication is fairly standard, and so we obtain the
following explicit version of the Erd\H{o}s-Stone-Bollob\'{a}s theorem.

\begin{theorem}
\label{ESnew}Let $r\geq2,$ let $c$ and $n$ be such that
\[
0<c<1\text{ \ \ and \ \ }n\geq\exp\left(  \left(  r^{r}/c\right)
^{r+1}\right)  ,
\]
and let $G$ be a graph of order $n.$ If $e\left(  G\right)  \geq\left(
1-1/r+c\right)  n^{2},$ then $G$ contains a $K_{r}\left(  s;t\right)  $ with
\[
s=\left\lfloor \left(  c/r^{r}\right)  ^{r+1}\log n\right\rfloor \text{
\ \ and \ \ }t>n^{1-\left(  c/r^{r}\right)  ^{r}}.
\]

\end{theorem}

In the two theorems above, we would like to emphasize the three principles
outlined in the introduction: first, the fundamental parameter $c$ may depend
on $n,$ e.g., letting $c=1/\log\log n,$ the conclusion is meaningful for
sufficiently large $n;$ note that this fact can be verified precisely because
the conditions for validity are stated explicitly. Also, the proof of these
theorems relies on more basic statements of wider applicability - Lemma
\ref{momo} and Lemma \ref{routl}.

Another observation about this setup is the peculiarity of the graphs
$K_{r}\left(  s;t\right)  $ in the conclusions of the above theorems: if the
statement holds for some $c,$ then it holds also for all positive $c^{\prime
}<c$ as long as $n$ is large enough. That is to say, when $n$ increases, in
addition to the graphs $K_{r}\left(  s;t\right)  $ guaranteed by the theorems,
we can find other, larger and more lopsided graphs $K_{r}\left(  s^{\prime
};t^{\prime}\right)  $ with $s^{\prime}<s$ and $t^{\prime}>t.$ This same
observation can be made on numerous other occasions below, and usually we
shall omit it to avoid repetition.

Let us note that Theorem \ref{BE} implies also the following assertion, which
strengthens the observation of Erd\H{o}s and Simonovits \cite{ErSi66}:

\begin{theorem}
\label{tESN} Let $r\geq3$ and let $F_{1},F_{2},\ldots$ be $\left(  r+1\right)
$-chromatic graphs satisfying $v\left(  F_{n}\right)  =o\left(  \log n\right)
.$ Then%
\[
\max\left\{  e\left(  G\right)  :G\in\mathcal{G}\left(  n\right)  \text{ and
}F_{n}\nsubseteq G\right\}  =\frac{r-1}{2r}n^{2}+o\left(  n^{2}\right)  .
\]

\end{theorem}

Thus Theorem \ref{tESN} solves asymptotically the Tur\'{a}n problem for
families of forbidden subgraphs whose order grows not too fast with $n.$
Moreover, the condition $v\left(  F_{n}\right)  =o\left(  \log n\right)  $ can
be sharpened further using the bounds given by Ishigami in \cite{Ish02}.

\subsubsection{\label{copy}Graphs with many copies of a given subgraph}

In this subsection we shall apply the basic idea above to arbitrary subgraphs
of graphs, including induced ones.

Let us first define a \emph{blow-up}\ of a graph $H$: given a graph $H$ of
order $r$ and positive integers $k_{1},\ldots,k_{r}$, we write $H\left(
k_{1},\ldots,k_{r}\right)  $ for the graph obtained by replacing each vertex
$u\in V\left(  H\right)  $ with a set $V_{u}$ of size $k_{u}$ and each edge
$uv\in E\left(  H\right)  $ with a complete bipartite graph with vertex
classes $V_{u}$ and $V_{v}.$

We are interested in the following generalization of Theorem \ref{Nik1}:
\emph{Suppose that a graph }$G$\emph{ of order }$n$\emph{ contains }$cn^{r}%
$\emph{ copies of a given subgraph }$H$\emph{ on }$r$\emph{ vertices. How
large a \textquotedblleft blow-up\textquotedblright\ of }$H$\emph{ must }%
$G$\emph{ contain?}\medskip

The following theorem from \cite{Nik08a} is an analog of Theorem \ref{Nik1}
for arbitrary subgraphs.

\begin{theorem}
\label{tsub}Let $r\geq2,$ let $c$ and $n$ be such that%
\[
0<c<1/r!\text{ \ \ and \ \ }n\geq\exp\left(  c^{r^{2}}\right)  ,
\]
and let $H$ be a graph of order $r.$ If $G\in\mathcal{G}\left(  n\right)  $
and $G$ contains more than $cn^{r}$ copies of $H,$ then $G$ contains an
$H\left(  s,\ldots s,t\right)  $ with $s=\left\lfloor c^{r^{2}}\log
n\right\rfloor $ and $t>n^{1-c^{r-1}}.$
\end{theorem}

A similar theorem is conceivable for induced subgraphs, but note the obvious
bump: the complete graph $K_{n}$ has $\Theta\left(  n^{2}\right)  $ edges,
i.e. $K_{2}$'s, but contains no \emph{induced }$4$-cycle, i.e. $K_{2}\left(
2\right)  .$ To come up with a meaningful statement, we need the following
more flexible version of a blow-up:\medskip

\textit{We say that a\ graph }$F$\textit{ \textbf{is of type} }$H\left(
k_{1},\ldots,k_{r}\right)  ,$\textit{ if }$F$\textit{ is obtained from
}$H\left(  k_{1},\ldots,k_{r}\right)  $\textit{ by adding some (possibly zero)
edges within the sets }$V_{u},$\textit{ }$u\in V\left(  H\right)  .$\medskip

This definition in hand, we can state the induced graph version of Theorem
\ref{tsub}, also from \cite{Nik08a}.

\begin{theorem}
Let $r\geq2,$ let $c$ and $n$ be such that%
\[
0<c<1/r!\text{ \ \ and \ \ }n\geq\exp\left(  c^{r^{2}}\right)  ,
\]
and let $H$ be a graph of order $r.$ If $G\in\mathcal{G}\left(  n\right)  $
and $G$ contains more than $cn^{r}$ induced copies of $H,$ then $G$ contains
an induced subgraph of type $H\left(  s,\ldots s,t\right)  ,$ where
$s=\left\lfloor c^{r^{2}}\log n\right\rfloor $ and $t>n^{1-c^{r-1}}.$
\end{theorem}

For constant $c$, the above theorems give the correct order of magnitude of
the subgraphs of type $H\left(  s,\ldots s,t\right)  ,$ namely, $\log n$ for
$s$ and $n^{1-o\left(  c\right)  }$ for $t.$ When $c$ depends on $n,$ the best
bounds on $s$ and $t$ are apparently unknown.

\subsubsection{\label{r-gr}Complete $r$-partite subgraphs of dense $r$-graphs}

In this subsection \emph{graph} stands for $r$\emph{-uniform hypergraph} for
some fixed $r\geq3$. We use again $K_{r}\left(  s_{1},\ldots,s_{r}\right)  $
to denote the complete $r$-partite $r$-graph with class sizes $s_{1}%
,\ldots,s_{r}.$

In the spirit of the previous topics, it is natural to ask: \emph{Suppose that
a graph }$G$\emph{ of order }$n$\emph{ contains }$cn^{r}$\emph{ edges. How
large a subgraph }$K_{r}\left(  s\right)  $ \emph{must} $G$ \emph{contain?} As
shown by Erd\H{o}s and Stone \cite{ErSt46} and Erd\H{o}s \cite{Erd64}, $s\geq
a\left(  \log n\right)  ^{1/\left(  r-1\right)  }$ for some $a=a\left(
c\right)  >0,$ independent of $n$.

In \cite{Nik09a} this fundamental result was extended in three directions: $c$
may depend on $n,$ the complete $r$-partite subgraph may have vertex classes
of variable size, and the graph $G$ is taken to be an $r$-partite $r$-graph
with equal classes. The last setup is obviously more general than just taking
$r$-graphs.

The following three theorems are given in \cite{Nik09a}.

\begin{theorem}
\label{tEr}Let $r\geq3,$ let $c$ and $n$ be such that
\[
0<c\leq r^{-3}\text{ \ \ and \ \ }n\geq\exp\left(  1/c^{r-1}\right)  ,
\]
and let the positive integers $s_{1},\ldots,s_{r-1}$ satisfy $s_{1}s_{2}\cdots
s_{r-1}\leq c^{r-1}\log n.$ Then every graph with $n$ vertices and at least
$cn^{r}/r!$ edges contains a $K_{r}\left(  s_{1},\ldots,s_{r-1},t\right)  $
with $t>n^{1-c^{r-2}}.$
\end{theorem}

Instead of this theorem it is easier and more effective to prove a more
general one for $r$-partite $r$-graphs.

\begin{theorem}
\label{tErd}Let $r\geq3,$ let $c$ and $n$ be such that
\[
0<c\leq r^{-3}\text{ \ \ and \ \ }n\geq\exp\left(  1/c^{r-1}\right)  ,
\]
and let the positive integers $s_{1},\ldots,s_{r-1}$ satisfy $s_{1}s_{2}\cdots
s_{r-1}\leq c^{r-1}\log n.$ Let $U_{1},\ldots,U_{r}$ be sets of size $n$ and
$E\subset U_{1}\times\cdots\times U_{r}$ satisfy $\left\vert E\right\vert \geq
cn^{r}.$ Then there exist $V_{1}\subset U_{1},\cdots,V_{r}\subset U_{r}$
satisfying $V_{1}\times\cdots\times V_{r}\subset E$ and \
\[
\left\vert V_{1}\right\vert =s_{1},\cdots,\left\vert V_{r-1}\right\vert
=s_{r-1},\text{ \ }\left\vert V_{r}\right\vert >n^{1-c^{r-2}}.
\]

\end{theorem}

In turn, Theorem \ref{tErd} is deduced from a counting result about
$r$-partite $r$-graphs, which generalizes the double counting argument of
K\"{o}vari, S\'{o}s and Tur\'{a}n for bipartite graphs \cite{KST54}.

\begin{theorem}
Let $r\geq2$ and let $c$ and $n$ be such that
\[
2^{r}\exp\left(  -\frac{1}{r}\left(  \log n\right)  ^{1/r}\right)  \leq
c\leq1.
\]
Let $G$ be an $r$-partite $r$-graph with parts $U_{1},\ldots,U_{r}$ of size
$n,$ and with edge set $E\subset U_{1}\times\cdots\times U_{r}$ satisfying
$\left\vert E\right\vert \geq cn^{r}.$ If the positive integers $s_{1}%
,s_{2},\ldots,s_{r}$ satisfy $s_{1}s_{2}\cdots s_{r}\leq\log n,$ then
$G\ $contains at least
\[
\left(  \frac{c}{2^{r}}\right)  ^{rs_{1}\cdots s_{r}}\binom{n}{s_{1}}%
\cdots\binom{n}{s_{r}}.
\]
complete $r$-partite subgraphs with precisely $s_{i}$ vertices in $U_{i}$ for
every $i=1,\ldots,r.$
\end{theorem}

Following Erd\H{o}s \cite{Erd64} and taking a random $r$-graph $G$ of order
$n$ and density $1-\varepsilon$, a straightforward calculation shows that with
probability tending to $1$, $G$ does not contain a $K_{r}\left(
s,\ldots,s\right)  $ for $s>A\left(  \log n\right)  ^{1/\left(  r-1\right)
},$ where $A=A\left(  \varepsilon\right)  $ is independent of $n.$ That is to
say, Theorems \ref{tEr} and \ref{tErd} are essentially tight.

\subsection{\label{sat}Saturation problems}

Saturation problems concern the type of subgraphs one necessarily finds in
graphs of order $n,$ with $t_{r}\left(  n\right)  +o\left(  n^{2}\right)  $
edges. Among all possible saturation problems we will consider only the most
important case: which subgraphs necessarily occur in graphs of order $n$ and
size $t_{r}\left(  n\right)  +1?$ Tur\'{a}n's theorem says that such graphs
contain a $K_{r+1},$ but one notes that they contain much larger supergraphs
of $K_{r+1}$.

Our first theorem completes an unfinished investigation started by Erd\H{o}s
in 1963, in \cite{Erd63}. We also present several results related to joints -
a class of important subgraphs, whose study was also initiated by Erd\H{o}s.

\subsubsection{Unavoidable subgraphs of graphs in $G\left(  n,t_{r}\left(
n\right)  +1\right)  $}

Let $s_{1}\geq2,$ and write $K_{r}^{+}\left(  s_{1},s_{2},...,s_{r}\right)  $
for the graph obtained from $K_{r}\left(  s_{1},s_{2},...,s_{r}\right)  $ by
adding an edge to the first part. For short, we also set
\[
K_{r}^{+}\left(  p\right)  =K_{r}^{+}\left(  p,...,p\right)  \text{ \ \ and
\ \ }K_{r}^{+}\left(  p;q\right)  =K_{r}^{+}\left(  p,...,p,q\right)  .
\]
In \cite{Erd63} Erd\H{o}s gave the following result:

\begin{theorem}
\label{Er2}For every $\varepsilon>0,$ there exist $c=c\left(  \varepsilon
\right)  >0$ and $n_{0}\left(  \varepsilon\right)  $ such that if $G$ is a
graph of order $n>n_{0}\left(  \varepsilon\right)  $ and $e\left(  G\right)
>\left\lfloor n^{2}/4\right\rfloor $, then $G$ contains a
\[
K_{2}^{+}\left(  \left\lfloor c\log n\right\rfloor ,\left\lceil
n^{1-\varepsilon}\right\rceil \right)  .
\]

\end{theorem}

For some time there was no generalization of this result for $K_{r}^{+}\left(
s;t\right)  $ until Erd\H{o}s and Simonovits \cite{ErSi73} came up with a
similar assertion valid for all $r\geq2$.

\begin{theorem}
\label{ESq}Let $r\geq2,$ $q\geq1,$ and let $n$ be sufficiently large. If $G$
is a graph of order $n$ with $t_{r}\left(  n\right)  +1$ edges, then $G$
contains a $K_{r}^{+}\left(  q\right)  .$
\end{theorem}

In a sense Theorem \ref{ESq} is best possible as any graph $H$ that
necessarily occurs in all sufficiently large graphs $G\in G\left(
n,t_{r}\left(  n\right)  +1\right)  $ can be imbedded in $K_{r}^{+}\left(
q\right)  $ for $q$ sufficiently large. To see this, just add an edge to the
Tur\'{a}n graph $T_{r}\left(  n\right)  $ and note that all $\left(
r+1\right)  $-partite subgraphs of this graph are edge-critical with respect
to the chromatic number. However, Theorem \ref{Er2} suggests that stronger
statements are possible, and indeed, in \cite{Nik10}, we extended both
Theorems \ref{Er2} and \ref{ESq} to the following one.

\begin{theorem}
\label{Th1}Let $r\geq2,$ let $c$ and $n$ be such that%
\[
0<c\leq r^{-\left(  r+7\right)  \left(  r+1\right)  }\text{ \ \ and \ \ }n\geq
e^{2/c},
\]
and let $G$ be a graph of order $n.$ If $e\left(  G\right)  >t_{r}\left(
n\right)  $, then $G$ contains a
\[
K_{r}^{+}\left(  \left\lfloor c\log n\right\rfloor ;\left\lceil n^{1-\sqrt{c}%
}\right\rceil \right)  .
\]

\end{theorem}

As usual, in Theorem \ref{Th1} $c$ may depend on $n$ within the given confine.
Note also that if the conclusion holds for some $c,$ it holds also for
positive $c^{\prime}<c,$ provided $n$ is sufficiently large. This implies
Erd\H{o}s's Theorem \ref{Er2}.

\subsubsection{Joints and books}

Erd\H{o}s \cite{Erd69} proved that if $r\geq2$ and $n>n_{0}\left(  r\right)
,$ every graph $G=G\left(  n,t_{r}\left(  n\right)  +1\right)  $ has an edge
that is contained in at least $n^{r-1}/\left(  10r\right)  ^{6r}$ cliques of
order $\left(  r+1\right)  .$ This fundamental fact seems so important, that
in \cite{BoNi04} we found it necessary to give the following
definition:\medskip

\textit{An }$r$\textit{-\textbf{joint} of size }$t$\textit{ is a collection of
}$t$\textit{ distinct }$r$\textit{-cliques sharing an edge.}\medskip

Note that two $r$-cliques of an $r$-joint may share up to $r-1$ vertices and
that for $r>3$ there may be many nonisomorphic $r$-joints of the same size. We
shall write $js_{r}\left(  G\right)  $ for the maximum size of an $r$-joint in
a graph $G$; in particular, if $2\leq r\leq n$ and $r$ divides $n,$ then
$js_{r}\left(  T_{r}(n)\right)  =\left(  \frac{n}{r}\right)  ^{r-2}$.

In this notation, the above result of Erd\H{o}s reads: \emph{if }$r\geq
2,$\emph{ }$n>n_{0}\left(  r\right)  ,$\emph{ and }$G\in\mathcal{G}\left(
n,t_{r}\left(  n\right)  +1\right)  ,$\emph{ then}%
\begin{equation}
js_{r+1}\left(  G\right)  \geq\frac{n^{r-1}}{\left(  10r\right)  ^{6r}}.
\label{Erdb}%
\end{equation}

In fact, the study of $js_{3}\left(  G\right)  $, also known as the
\emph{booksize} of $G$, was initiated by Erd\H{o}s even earlier, in
\cite{Erd62b}, and was subsequently generalized in \cite{Erd68} and
\cite{Erd69}; it seems that he foresaw the importance of joints when he
restated his general result in 1995, in \cite{EFGG95}. A quintessential result
concerning joints is the \textquotedblleft triangle removal
lemma\textquotedblright\ of Ruzsa and Szemer\'{e}di \cite{RuSz75}, which can
be stated as a lower bound on the booksize $js_{3}\left(  G\right)  $ when $G$
is a graph of a particular kind.

In fact joints help to obtain several of the results mentioned in this survey,
e.g., the general stability Theorem \ref{tStab} and its spectral version,
Theorem \ref{stStab}. Later, we shall give also spectral conditions for the
existence of large joints, in Theorem \ref{spN1}.

In \cite{BoNi04}, Bollob\'{a}s and the author enhanced the bound of Erd\H{o}s
(\ref{Erdb}) to the following explicit one.

\begin{theorem}
\label{tJ} Let $r\geq2,$ $n>r^{8},$ and let $G$ be a graph of order $n.$ If
$e\left(  G\right)  \geq t_{r}\left(  n\right)  ,$ then%
\[
js_{r+1}\left(  G\right)  >\frac{n^{r-1}}{r^{r+5}}%
\]
unless $G=T_{r}\left(  n\right)  .$
\end{theorem}

In \cite{BoNi10} an analogous theorem is given in the case when $G$ has many
$r$-cliques, rather than edges. More precisely, letting $k_{r}\left(
G\right)  $ stand for the number of $r$-cliques of a graph $G,$ we have

\begin{theorem}
\label{tJr} Let $r\geq2,$ $n>r^{8},$ and let $G$ be a graph of order $n.$ If
$k_{s}(G)\geq k_{s}\left(  T_{r}(n)\right)  $ for some $s,$ $\left(  2\leq
s\leq r\right)  ,$ then
\[
js_{r+1}\left(  G\right)  >\frac{n^{r-1}}{r^{2r+12}}%
\]
unless $G=T_{r}\left(  n\right)  $.
\end{theorem}

Note that Theorems \ref{tJ} and \ref{tJr} cannot be improved too much, as
shown by the graph $G$ obtained by adding an edge to $T_{r}(n)$: we have
$k_{s}(G)\geq k_{s}(T_{r}(n))$ but $js_{r+1}(G)\leq\left\lceil n/r\right\rceil
^{r-1}$. However, the best bound in Theorem \ref{tJ} is known only for
$3$-joints. Usually a $3$-joint of size $t$ is called a \emph{book of size}
$t.$ Edwards \cite{EdMS}, and independently Khad\v{z}iivanov and Nikiforov
\cite{KhNi79} proved the following theorem.

\begin{theorem}
If $G$ is a graph of order $n$ with $e\left(  G\right)  >\left\lfloor
n^{2}/4\right\rfloor ,$ then it contains a book of size greater than $n/6.$
\end{theorem}

This theorem is best possible in view of the following graph. Let $n=6k.$
Partition $\left[  n\right]  $ into $6$ sets $A_{11},$ $A_{12},$ $A_{13},$
$A_{21},$ $A_{22},$ $A_{23}$ with $\left\vert A_{11}\right\vert =\left\vert
A_{12}\right\vert =\left\vert A_{13}\right\vert =k-1$ and $\left\vert
A_{21}\right\vert =\left\vert A_{22}\right\vert =\left\vert A_{23}\right\vert
=k+1.$ For $1\leq j<k\leq3$ join every vertex of $A_{ij}$ to every vertex of
$A_{ik}$ and for $j=1,2,3$ join every vertex of $A_{1j}$ to every vertex of
$A_{2j}.$ The resulting graph has size $>\left\lfloor n^{2}/4\right\rfloor +1$
and its booksize is $k+1=n/6+1.$

A more recent presentation of these results can be found in \cite{BoNi05}.

\subsection{\label{stab}Stability problems}

This subsection has three parts. First we sharpen the classical stability
theorem of Erd\H{o}s \cite{Erd66},\cite{Erd68} and Simonovits \cite{Sim68},
which gives information about the structure of graphs without fixed forbidden
subgraphs and whose size is close to the maximum possible. Second, we give
several specific stability theorems for specific forbidden subgraphs, where
stronger conclusions are possible. Lacking a better term, we call such cases
\emph{strong stability}.

Finally, we discuss the structure of $K_{r}$-free graphs of large minimum
degree. This is a rich area with many results and a long history. It is not
customary to consider it in the context of stability problems, but we believe
this is the general category where this area belongs, since most of its
statements can be phrased so that large minimum degree of a $K_{r}$-free graph
implies a certain structure.

\subsubsection{A general stability theorem}

Let $F$ be a fixed $\left(  r+1\right)  $-partite graph $F$ and $G$ be a graph
of order $n.$ The theorem of Erd\H{o}s and Stone implies that if
$\varepsilon>0$ and $e\left(  G\right)  >\left(  1-1/r+\varepsilon\right)
n^{2}/2,$ then $G$ contains $F,$ when $n$ is sufficiently large. On the other
hand, $T_{r}\left(  n\right)  $ is $r$-partite and therefore does not contain
$F,$ although
\[
e\left(  T_{r}\left(  n\right)  \right)  =t_{r}\left(  n\right)
\thickapprox\left(  1-1/r\right)  n^{2}/2.
\]
Erd\H{o}s and Simonovits \cite{Erd66},\cite{Erd68},\cite{Sim68} noticed that
if a graph $G$ of order $n$ contains no copy of $F$ and has close to $\left(
1-1/r\right)  n^{2}/2$ edges, then $G$ is similar to $T_{r}\left(  n\right)  $.

\begin{theorem}
\label{ESstab}Let $r\geq2$ and let $F$ be a fixed $\left(  r+1\right)
$-partite graph. For every $\delta>0,$ there is an $\varepsilon>0$ such that
if $G$ is a graph of order $n$ with $e\left(  G\right)  >\left(
1-1/r-\varepsilon\right)  n^{2}/2,$ then either $G$ contains $F$ or $G$
differs from $T_{r}\left(  n\right)  $ in fewer than $\delta n^{2}$ edges.
\end{theorem}

A closer inspection of this statement reveals that $\varepsilon$ depends both
on $\delta$ and on $F.$ To investigate this dependence, we simplify the
picture by assuming that $F$ is a complete $\left(  r+1\right)  $-graph.
Moreover, radically departing from the setup of fixed $F$, we assume that
$F=K_{r+1}\left(  \left\lfloor c\log n\right\rfloor ;\left\lceil n^{1-\sqrt
{c}}\right\rceil \right)  $ for some $c>0.$ Note that for a given $n$ the
single real parameter $c$ characterizes $F$ completely. It turns out with this
selection of $F$ we still can get an enhancement of Theorem \ref{ESstab}, as
proved in \cite{Nik09c}.

\begin{theorem}
\label{tStab}Let $r\geq2,$ let $c,$ $\varepsilon$ and $n$ be such that
\[
0<c<r^{-3\left(  r+14\right)  \left(  r+1\right)  },\text{ \ \ }%
0<\varepsilon<r^{-24},\text{ \ \ }n>e^{1/c},
\]
and let $G$ be a graph of order $n.$ If $e\left(  G\right)  >\left(
1-1/r-\varepsilon\right)  n^{2}/2$, then one of the following statements holds:

(a) $G$ contains a $K_{r+1}\left(  \left\lfloor c\log n\right\rfloor
;\left\lceil n^{1-\sqrt{c}}\right\rceil \right)  ;$

(b) $G$ differs from $T_{r}\left(  n\right)  $ in fewer than $\left(
\varepsilon^{1/3}+c^{1/\left(  3r+3\right)  }\right)  n^{2}$ edges.
\end{theorem}

Note that, as usual, $c$ may depend on $n.$ A natural question is how tight
Theorem \ref{tStab} is. The complete answer seems difficult since two
parameters, $\varepsilon$ and $c,$ are involved. First, the factor $\left(
\varepsilon^{1/3}+c^{1/\left(  3r+3\right)  }\right)  $ in condition
\emph{(b)} is far from the best one, but is simple. However for fixed $c$
condition \emph{(a) }is best possible up to a constant factor. Indeed, let
$\alpha>0$ be sufficiently small. A randomly chosen graph of order $n$ with
$\left(  1-\alpha\right)  n^{2}/2$ edges contains no $K_{2}\left(
\left\lfloor c^{\prime}\log n\right\rfloor ,\left\lfloor c^{\prime}\log
n\right\rfloor \right)  $ and differs from $T_{r}\left(  n\right)  $ in more
that $c^{\prime\prime}n^{2}$ edges for some positive $c^{\prime}$ and
$c^{\prime\prime},$ independent of $n$.

\subsubsection{Strong stability}

For certain forbidden graphs condition \emph{(ii)} of Theorem \ref{tStab} can
be strengthened. Such particular stability theorems can be of interest in
applications, e.g., Ramsey problems. We start with a theorem in \cite{NiRo05},
which gives a particular stability condition for $K_{r+1}$-free graphs.

\begin{theorem}
\label{stabK}Let $r\geq2$ and $0<\varepsilon\leq2^{-10}r^{-6}$, and let $G$ be
a $K_{r+1}$-free graph of order $n.$ If $e\left(  G\right)  >\left(
1-1/r-\varepsilon\right)  n^{2}/2,$ then $G$ contains an induced $r$-partite
graph $H$ of order at least $\left(  1-2\sqrt[3]{\varepsilon}\right)  n$ and
with minimum degree $\delta\left(  H\right)  \geq\left(  1-1/r-4\sqrt[3]%
{\varepsilon}\right)  n.$
\end{theorem}

Note that the stability condition in this theorem is stronger than condition
\emph{(b)} of Theorem \ref{tStab}. Indeed, the classes of $H$ are almost
equal, it is almost complete, and contains almost all vertices of $G.$ This
type of conclusion is the purpose of the three theorems below. In the first
two of them the premise \textquotedblleft$K_{r+1}$-free\textquotedblright%
\ will be further weakened; but Theorem \ref{stabK} is still of interest,
because it is proved for all conceivable $n$.

The following two theorems have been proved in \cite{Nik10} and \cite{BoNi04}.

\begin{theorem}
\label{Th2}Let $r\geq2,$ let $c,$ $\varepsilon$ and $n$ be such that
\[
0<c<r^{-\left(  r+7\right)  \left(  r+1\right)  }/2,\text{ \ \ }%
0<\varepsilon<r^{-8}/8,\text{ \ \ }n>e^{2/c},
\]
and let $G$ be a graph of order $n.$ If $e\left(  G\right)  >\left(
1-1/r-\varepsilon\right)  n^{2}/2$, then one of the following statements holds:

(a) $G$ contains a $K_{r}^{+}\left(  \left\lfloor c\log n\right\rfloor
;\left\lceil n^{1-2\sqrt{c}}\right\rceil \right)  ;$

(b) $G$ contains an induced $r$-partite subgraph $H$ of order at least
$\left(  1-\sqrt{2\varepsilon}\right)  n,$ with minimum degree
\[
\delta\left(  H\right)  >\left(  1-1/r-2\sqrt{2\varepsilon}\right)  n.
\]

\end{theorem}

\begin{theorem}
\label{stabj} Let $r\geq2,$ let $c$ and $n$ be such that
\[
r\geq2,\text{ \ \ }0<\varepsilon<r^{-8}/32,\text{ \ \ }n>r^{8},
\]
and let $G$ be a graph of order $n.$ If $e\left(  G\right)  >\left(
1-1/r-\varepsilon\right)  n^{2}/2,$ then one of the following statements holds:

(a) $js_{r+1}\left(  G\right)  >\left(  1-1/r^{3}\right)  n^{r-1}/r^{r+5};$

(b) $G$ contains an induced $r$-partite subgraph $H$ of order at least
$\left(  1-4\sqrt{\varepsilon}\right)  n,$ with minimum degree
\[
\delta\left(  H\right)  >\left(  1-1/r-6\sqrt{\varepsilon}\right)  n.
\]

\end{theorem}

As one can expect, the analogous statement for books is quite close to the
best possible \cite{BoNi05}.

\begin{theorem}
\label{stabb} Let $0<\varepsilon<10^{-5}$ and let $G$ be a graph of order $n.$
If $e\left(  G\right)  >\left(  1/4-\varepsilon\right)  n^{2},$ then either
$G$ contains a book of size at least $\left(  1/6-2\sqrt[3]{\varepsilon
}\right)  n$ or $G$ contains an induced bipartite graph $H$ of order at least
$\left(  1-\sqrt[3]{\varepsilon}\right)  n$ and with minimal degree
$\delta\left(  H\right)  \geq\left(  1/2-4\sqrt[3]{\varepsilon}\right)  n.$
\end{theorem}

\subsubsection{\label{3free}$K_{r}$-free graphs with large minimum degree}

A famous theorem of Andr\'{a}sfai, Erd\H{o}s and S\'{o}s \cite{AES74} shows
that if $r\geq2$ and $G$ is a $K_{r+1}$-free graph of order $n$ and with
minimum degree satisfying%
\begin{equation}
\delta\left(  G\right)  >\left(  1-\frac{3}{3r-1}\right)  n, \label{AESbas}%
\end{equation}
then $G$ is $r$-partite. They also gave an example showing that equality in
(\ref{AESbas}) is not sufficient to get the same conclusion.

In particular, for $r=2$ this statement says that every triangle-free graph of
order $n$ with minimum degree $\delta\left(  G\right)  >2n/5$ is bipartite. On
the other hand, Hajnal \cite{ErSi73} constructed a triangle-free graph of
order $n$ with arbitrary large chromatic number and with minimum degree
$\delta\left(  G\right)  >\left(  1/3-\varepsilon\right)  n.$ In view of
Hajnal's example, Erd\H{o}s and Simonovits \cite{ErSi73} conjectured that all
$K_{3}$-free graphs of order $n$ with $\delta\left(  G\right)  >n/3$ are
$3$-chromatic. However, this conjecture was disproved by H\"{a}ggkvist
\cite{Hag82}, who described for every $k\geq1$ a $10k$-regular, $4$-chromatic,
triangle-free graph of order $29k.$ The example of H\"{a}ggkvist is based on
the Mycielski graph $M_{3},$ also known as the Gr\"{o}tzsch graph, which is a
$4$-chromatic triangle-free graph of order $11$. To construct $M_{3},$ let
$v_{1},\ldots,v_{5}$ be the vertices of a $5$-cycle and choose $6$ other
vertices $u_{1},\ldots,u_{6}.$ Join $u_{i}$ to the neighbors of $v_{i}$ for
all $i=1,\ldots,5$, and finally join $u_{6}$ to $u_{1},\ldots,u_{5}.$

Other graphs that are crucial in these questions are the triangle-free,
$3$-chromatic Andr\'{a}sfai graphs $A_{1},A_{2},\ldots,$ first described in
\cite{And62}:\ \emph{set }$A_{1}=K_{2}$\emph{ and for every }$i\geq2$\emph{
let }$A_{i}$\emph{ be the complement of the }$\left(  i-1\right)  $\emph{-th
power of the cycle }$C_{3i-1}.$

To state the next structural theorems we need the following definition:
\emph{a graph }$G$\emph{ is said to be \textbf{homomorphic} to a graph }%
$H,$\emph{ if there exists a map }$f:V\left(  G\right)  \rightarrow V\left(
H\right)  $\emph{ such that }$uv\in E\left(  G\right)  $\emph{ implies that
}$f\left(  u\right)  f\left(  v\right)  \in E\left(  H\right)  $\emph{.}

In \cite{Jin93}, Jin generalized the case $r=2$ of the theorem of
Andr\'{a}sfai, Erd\H{o}s and S\'{o}s and a result of H\"{a}ggkvist from
\cite{Hag82} in the following theorem.

\begin{theorem}
\label{tA}Let $1\leq k\leq9,$ and let $G$ be a triangle-free graph of order
$n.\ $If
\[
\delta\left(  G\right)  >\frac{k+1}{3k+2}n,
\]
then $G$ is homomorphic to $A_{k}.$
\end{theorem}

Note that this result is tight: taking the graph $A_{k+1}$, and blowing it up
by a factor $t,$ we obtain a triangle-free graph $G$ of order $n=\left(
3k+2\right)  t$ vertices, with $\delta\left(  G\right)  =\left(  k+1\right)
n/\left(  3k+2\right)  ,$ which is not homomorphic to $A_{k}.$

Note also that all graphs satisfying the premises of Theorem\ \ref{tA} are
$3$-chromatic. Addressing this last issue, Jin \cite{Jin95}, and Chen, Jin and
Koh \cite{CJK97} gave a finer characterization of all $K_{3}$-free graphs
with\emph{ }$\delta>n/3.$

\begin{theorem}
\label{tB}Let $G$ be a triangle-free graph of order $n,$ with $\delta\left(
G\right)  >n/3.$ If $\chi\left(  G\right)  \geq4,$ then $M_{3}\subset G.$ If
$\chi\left(  G\right)  =3$ and
\[
\delta\left(  G\right)  >\frac{k+1}{3k+2}n,
\]
then $G$ is homomorphic to $A_{k}.$
\end{theorem}

Finally, Brandt and Thomass\'{e} \cite{BrTo10} gave the following ultimate result.

\begin{theorem}
\label{tC} Let $G$ be a triangle-free graph of order $n.$ If $\delta\left(
G\right)  >n/3,$ then $\chi\left(  G\right)  \leq4.$
\end{theorem}

It is natural to ask the same questions for $K_{r}$-free graphs with large
minimum degree. Contrary to expectation, the answers are by far easier. First,
the graphs of Andr\'{a}sfai, Hajnal and H\"{a}ggkvist are easily generalized
by joining them with appropriately chosen complete $\left(  r-3\right)
$-partite graphs.

In particular, for every $\varepsilon$ there exists a $K_{r+1}$-free graph of
order $n$ with
\[
\delta\left(  G\right)  >\left(  1-\frac{2}{2r-1}-\varepsilon\right)  n
\]
and arbitrary large chromatic number, provided $n$ is sufficiently large.

Hence, the main question is: \emph{how large }$\chi\left(  G\right)  $
\emph{can be when }$G$\emph{ is a }$K_{r+1}$\emph{-free graph of order }%
$n$\emph{ }$\emph{with}$ $\delta\left(  G\right)  >\left(  1-2/\left(
2r-1\right)  \right)  n.$ The answer is:

\begin{theorem}
\label{extBT} Let $r\geq2$ and let $G$ be a $K_{r+1}$-free graph of order $n.$
If
\[
\delta\left(  G\right)  >\left(  1-\frac{2}{2r-1}\right)  n,
\]
then $\chi\left(  G\right)  \leq r+2.$
\end{theorem}

This theorem leaves only two cases to investigate, viz., $\chi\left(
G\right)  =r+1$ and $\chi\left(  G\right)  =r+2$. As one can expect, when
$\delta\left(  G\right)  $ is sufficiently large, we have $\chi\left(
G\right)  =r+1.$ The precise statement extends Theorem \ref{tA} as follows.

\begin{theorem}
\label{extJin}Let $r\geq2,$ $1\leq k\leq9,$ and let $G$ be a $K_{r+1}$-free
graph of order $n.$ If
\[
\delta\left(  G\right)  >\left(  1-\frac{2k-1}{\left(  2k-1\right)
r-k+1}\right)  n
\]
then $G$ is homomorphic to $A_{k}+K_{r-2}.$
\end{theorem}

As a corollary, under the premises of Theorem \ref{extJin}, we find that
$\chi\left(  G\right)  \leq r+1.$ Also Theorem \ref{extJin} is best possible
in the following sense: for every $k$ and $n,$ there exists an $\left(
r+1\right)  $-chromatic $K_{r+1}$-free $G$ of order $n$ with%
\[
\delta\left(  G\right)  \geq\left(  1-\frac{2k-1}{\left(  2k-1\right)
r-k+1}\right)  n-1
\]
that is not homomorphic to $A_{k}+K_{r-2}$.

Using the example of H\"{a}ggkvist, we construct for every $n$ an $\left(
r+2\right)  $-chromatic, $K_{r+1}$-free graph $G$ with
\[
\delta\left(  G\right)  \geq\left(  1-\frac{19}{19r-9}\right)  n-1,
\]
which shows that the conclusion of Theorem \ref{extJin} does not necessarily
hold for $k\geq10.$

To give some further structural information, we extend Theorem \ref{tC} as follows.

\begin{theorem}
\label{extCJK}Let $r\geq2$ and $G$ be a $K_{r+1}$-free graph of order $n$
with
\[
\delta\left(  G\right)  >\left(  1-\frac{2}{2r-1}\right)  n.
\]
If $\chi\left(  G\right)  \geq r+2,$ then $M_{3}+K_{r-2}\subset G.$ If
$\chi\left(  G\right)  \leq r+1$ and%
\[
\delta\left(  G\right)  >\left(  1-\frac{2k-1}{\left(  2k-1\right)
r-k+1}\right)  n
\]
then $G$ is homomorphic to $A_{k}+K_{r-2}.$
\end{theorem}

This result is best possible in view of the examples described prior to
Theorem \ref{extCJK}.

We deduce the proofs of Theorems \ref{extBT}, \ref{extJin} and \ref{extCJK} by
induction on $r$ from Theorems \ref{tC}, \ref{tA} and \ref{tB} respectively.
The induction step, carried out uniformly in all the three proofs, is based on
the crucial Lemma \ref{reduL}. This lemma can be applied immediately to extend
other results about triangle-free graphs.

The new results in this subsection, together with Lemma \ref{reduL} have been
published in \cite{Nik10c}. Since the first version of that paper was made
public, the author learned that similar research has been undertaken
independently by W. Goddard and J. Lyle \cite{GoLy10}.

\subsection{The number of graphs with large forbidden subgraphs}

An intriguing question is how many graphs with given properties there are. For
certain natural properties such as \textquotedblleft$G$ is $K_{r}%
$-free\textquotedblright\ or \textquotedblleft$G$ has no induced graph
isomorphic to $H$\textquotedblright\ satisfactory answers have been obtained.
Thus, given a graph $H,$ write $\mathcal{P}_{n}\left(  H\right)  $ for the set
of all labelled graphs of order $n$ not containing $H.$ A classical result of
Erd\H{o}s, Kleitman and Rothschild \cite{EKR73} states that
\begin{equation}
\log_{2}\left\vert \mathcal{P}_{n}\left(  K_{r+1}\right)  \right\vert =\left(
1-1/r+o\left(  1\right)  \right)  \binom{n}{2}. \label{EKR}%
\end{equation}
Ten years later, Erd\H{o}s, Frankl and R\"{o}dl \cite{EFR86} showed that the
conclusion in (\ref{EKR}) remains valid if $K_{r+1}$ is replaced by an
arbitrary fixed $\left(  r+1\right)  $-chromatic graph $H$.

In fact, as shown in \cite{BoNi09b}, the conclusion in (\ref{EKR}) also
remains valid if $K_{r+1}$ is replaced by a sequence of forbidden graphs whose
order grows with $n$. Until recently such results seemed to be out of reach;
however, the framework laid out in \cite{Nik08} and \cite{Nik08a} has opened
new possibilities. Here is the theorem that directly generalizes the
Erd\H{o}s-Frankl-R\"{o}dl result.

\begin{theorem}
\label{BN}Given $r\geq2$ and $0<\varepsilon\leq1/2,$ there exists
$\delta=\delta\left(  \varepsilon\right)  >0$ such that for $n$ sufficiently
large,%
\begin{equation}
\left(  1-1/r\right)  \binom{n}{2}\leq\log_{2}\left\vert \mathcal{P}%
_{n}\left(  K_{r+1}\left(  \left\lfloor \delta\log n\right\rfloor ;\left\lceil
n^{1-\sqrt{\delta}}\right\rceil \right)  \right)  \right\vert \leq\left(
1-1/r+\varepsilon\right)  \binom{n}{2}. \label{mi1}%
\end{equation}

\end{theorem}

Note that the real contribution of Theorem \ref{BN} is the upper bound in
(\ref{mi1}) since the lower bound follows by counting the labelled spanning
subgraphs of the Tur\'{a}n graph $T_{r}\left(  n\right)  .$ Let us mention
that the proof of Theorem \ref{BN} does not use Szemer\'{e}di's Regularity
Lemma, which is a standard tool for such questions.

Similar statements can be proved for forbidden induced subgraphs, where the
role of the chromatic number is played by the \emph{coloring number }$\chi
_{c}$ of a graph property, introduced first in \cite{BoTh95}, and defined
below.\medskip

\textit{Let }$0\leq s\leq r$\textit{ be integers and let }$\mathcal{H}\left(
r,s\right)  $\textit{ be the class of graphs whose vertex sets can be
partitioned into }$s$\textit{ cliques and }$r-s$\textit{ independent sets.
Given a graph property }$\mathcal{P}$\textit{, the coloring number }$\chi
_{c}\left(  \mathcal{P}\right)  $\textit{ is defined as}%
\[
\chi_{c}\left(  H\right)  =\max\left\{  r:\mathcal{H}\left(  r,s\right)
\subset\mathcal{P}\text{ \textit{for some }}s\in\left[  r\right]  \right\}
\]

Also, given a graph $H,$ let us write $\mathcal{P}_{n}^{\ast}\left(  H\right)
$ for the set of graphs of order $n$ not containing $H$ as an induced
subgraph. Alexeev \cite{Ale93} and, independently Bollob\'{a}s and Thomason
\cite{BoTh95},\cite{BoTh97} proved that the exact analog of the result of
Erd\H{o}s, Frankl and R\"{o}dl holds:\medskip

\emph{ If }$H$\emph{ is a fixed graph and }$r=\chi_{c}\left(  \mathcal{P}%
_{n}^{\ast}\left(  H\right)  \right)  $\emph{, then}
\begin{equation}
\log_{2}\left\vert \mathcal{P}_{n}^{\ast}\left(  H\right)  \right\vert
=\left(  1-1/r+o\left(  1\right)  \right)  \binom{n}{2}. \label{ABT}%
\end{equation}

This result also can be extended by replacing $H$ with a sequence of forbidden
graphs whose order grows with $n.$ To this end, recall the definition of a
graph \emph{of type }$H\left(  k_{1},\ldots,k_{h}\right)  ,$ where $H$ is a
fixed labelled graph of order $h$ and $k_{1},\ldots,k_{h}$ are positive
integers (Subsection \ref{copy})$.$ Informally, a graph of type $H\left(
k_{1},\ldots,k_{h}\right)  $ is obtained by first \textquotedblleft
blowing-up\textquotedblright\ $H$ to $H\left(  k_{1},\ldots,k_{h}\right)  $
and then adding (possibly zero) edges to the vertex classes of the
\textquotedblleft blow-up\textquotedblright\ but keeping intact the edges
across vertex classes.

Now, given a labelled graph $H$ and positive integers $p$ and $q,$ let
$\mathcal{P}_{n}\left(  H;p,q\right)  $ be the set of labelled graphs of order
$n$ that contain no induced subgraph of type $H\left(  p,\ldots,p,q\right)  $.

Here is the result for forbidden induced subgraphs, also from \cite{BoNi09b}.

\begin{theorem}
\label{BN2}Let $H$ be a labelled graph and let $r=\chi_{c}\left(
\mathcal{P}_{n}^{\ast}\left(  H\right)  \right)  .$ For every $\varepsilon>0,$
there is $\delta=\delta\left(  \varepsilon\right)  >0$ such that for $n$
sufficiently large%
\begin{equation}
\left(  1-1/r\right)  \binom{n}{2}\leq\log_{2}\left\vert \mathcal{P}%
_{n}\left(  H;\left\lfloor \delta\log n\right\rfloor ,\left\lceil
n^{1-\sqrt{\delta}}\right\rceil \right)  \right\vert \leq\left(
1-1/r+\varepsilon\right)  \binom{n}{2}. \label{bd1}%
\end{equation}

\end{theorem}

In a sense Theorems \ref{BN} and \ref{BN2} are almost best possible, in view
of the following simple observation, which can be proved by considering the
random graph $G_{n,p}$ with $p\rightarrow1.$\medskip

\emph{Given }$r\geq2$\emph{ and }$\varepsilon>0,$\emph{ there exists }%
$C>0$\emph{ such that the number }$S_{n}$ \emph{of labelled graphs which do
not contain }$K_{2}\left(  \left\lceil C\log n\right\rceil ,\left\lceil C\log
n\right\rceil \right)  $\emph{ satisfies }$S_{n}\geq\left(  1-\varepsilon
\right)  2^{\binom{n}{2}}.$

\section{\label{spec}Spectral extremal graph theory}

Generally speaking, spectral graph theory investigates graphs using the
spectra of various matrices associated with graphs, such as the adjacency
matrix. For an introduction to this topic we refer the reader to \cite{CRS10}.

Given a graph $G$ with vertex set $\left\{  v_{1},\ldots,v_{n}\right\}  ,$ the
adjacency matrix of $G$ is a matrix $A=\left[  a_{ij}\right]  $ of size $n$
given by%
\[
a_{ij}=\left\{
\begin{array}
[c]{ll}%
1, & \text{if }\left(  v_{i},v_{j}\right)  \in E\left(  G\right)  ;\\
0, & \text{otherwise.}%
\end{array}
\right.
\]

Note that $A$ is symmetric and nonnegative, and much is known about the
spectra of such matrices. For instance, the eigenvalues of $A$ are real
numbers, which we shall denote by $\mu_{1}\left(  G\right)  ,\ldots,\mu
_{n}\left(  G\right)  ,$ indexed in non-increasing order. The value
$\mu\left(  G\right)  =\mu_{1}\left(  G\right)  $ is called the spectral
radius of $G$ and has maximum absolute value among all eigenvalues.

Another matrix that we shall use is the Laplacian matrix $L,$ defined as
$D\left(  G\right)  -A\left(  G\right)  ,$ where $D\left(  G\right)  $ is the
diagonal matrix of the row-sums of $A,$ i.e., the degrees of $G.$ The
eigenvalues of the Laplacian are denoted by $\lambda_{1}\left(  G\right)
,\ldots,\lambda_{n}\left(  G\right)  ,$ indexed in non-decreasing order.

A third matrix associated with graphs is the $Q$-matrix or the
\textquotedblleft signless Laplacian\textquotedblright, defined as $Q=D+A.$
The eigenvalues of $Q$ are denoted by $q_{1}\left(  G\right)  ,\ldots
,q_{n}\left(  G\right)  ,$ indexed in non-increasing order. The $Q$-matrix has
received a lot of attention in recent years, see, e.g., \cite{CvSi09}%
,\cite{CvSi10a} and \cite{CvSi10b}. The Laplacian and the $Q$-matrix are
positive semi-definite matrices, and $\lambda_{1}\left(  G\right)  =0.$

\subsection{The spectral problems that are studied}

How large can be the spectral radius $\mu\left(  G\right)  $ when $G$ is a
$K_{r}$-free graph of order $n$? Such questions come easily to the mind when
one studies extremal graph problems. In fact, with any extremal problem of the
type \textquotedblleft\emph{What is the maximum number of edges in a graph
}$G$\emph{ of order }$n$ \emph{with property }$\mathcal{P}$\emph{?}%
\textquotedblright\emph{ }goes a spectral analog: \textquotedblleft\emph{What
is the maximum spectral radius} \emph{of a graph }$G$\emph{ of order }$n$
\emph{with property }$\mathcal{P}$\emph{?}\textquotedblright\ This is not
merely a superficial analogy since if we have a solution of the spectral
problem, then by the fundamental inequality
\begin{equation}
\mu\left(  G\right)  \geq2e\left(  G\right)  /n, \label{cosi}%
\end{equation}
we immediately obtain an upper bound on $e\left(  G\right)  $ as well. The use
of this implication is illustrated on several occasions below; in particular,
for the Zarankiewicz problem we obtain the sharpest bounds on $e\left(
G\right)  $ known so far.

On the other hand, inequality (\ref{cosi}) suggests a way to conjecture
spectral results by taking known nonspectral extremal statements that involve
the average degree of a graph and replacing the average degree by $\mu\left(
G\right)  .$ More often than not, the resulting statement is correct and even
stronger, but of course it needs its own proof. To create suitable proof tools
we painstakingly built several technical but rather flexible statements such
as Theorem \ref{thv4} and Lemma \ref{lev1}.

This smooth machinery is sufficient to prove spectral analogs of most of the
extremal problems discussed in Section \ref{class} and of several others as
well. Among these results are: various forms of Tur\'{a}n's theorem, the
Erd\H{o}s-Stone-Bollob\'{a}s theorem, conditions for large joints and for odd
cycles; a general stability theorem and several strong stability theorems, an
asymptotic solution of the general extremal problem for nonbipartite forbidden
subgraphs, the Zarankiewicz problem, sufficient conditions for paths and
cycles, sufficient conditions for Hamilton paths and cycles.

Despite these successful translations, more can be expected from spectral
extremal graph theory, which seems inherently richer than the conventional
one. Indeed, we give also a fair number of spectral results that have no
conventional analog, for example, results involving the smallest eigenvalue of
the adjacency matrix or the spectral radius of the Laplacian matrix.

\subsection{Spectral forms of the Tur\'{a}n theorem}

In 1986, Wilf \cite{Wil86} showed that if $G$ is a graph of order $n$ with
clique number $\omega\left(  G\right)  =\omega,$ then
\begin{equation}
\mu\left(  G\right)  \leq\left(  1-1/\omega\right)  n. \label{bnd4}%
\end{equation}

Note first that in view of the inequality $\mu\left(  G\right)  \geq2e\left(
G\right)  /n,$ (\ref{bnd4}) implies the concise Tur\'{a}n theorem:%
\begin{equation}
e\left(  G\right)  \leq\left(  1-1/\omega\right)  n^{2}/2. \label{conTT}%
\end{equation}
However, inequality (\ref{bnd4}) opens many other new possibilities. Indeed,
if we combine (\ref{bnd4}) with other lower bounds on $\mu(G)$, e.g., with
\[
\mu^{2}(G)\geq\frac{1}{n}\sum_{u\in V(G)}d^{2}\left(  u\right)  ,
\]
we obtain other forms of (\ref{conTT}). An infinite class of similar lower
bounds is given in \cite{Nik06}.

Below we sharpen inequality (\ref{bnd4}) in two ways.

\subsubsection*{A concise spectral Tur\'{a}n theorem}

In 1970 Nosal \cite{Nos70} showed that every triangle-free graph $G$ satisfies
$\mu\left(  G\right)  \leq\sqrt{e\left(  G\right)  }.$ This result was
extended in \cite{Nik02} and \cite{Nik09d} in the following theorem,
conjectured by Edwards and Elphick in \cite{EdEl83}:

\begin{theorem}
\label{Nik02}If $G$ is a graph of order $n$ and $\omega\left(  G\right)
=\omega,$ then%
\begin{equation}
\mu^{2}\left(  G\right)  \leq2\left(  1-1/\omega\right)  e\left(  G\right)  .
\label{bnd6}%
\end{equation}
If $G$ has no isolated vertices, then equality is possible if and only if one
of the following conditions holds:

(a) $\omega=2$ and $G$ is a complete bipartite graph;

(b) $\omega\geq3$ and $G$ is a complete regular $\omega$-partite graph.
\end{theorem}

In view of (\ref{conTT}), we see that%
\[
\mu^{2}\left(  G\right)  \leq2\left(  1-1/\omega\right)  m\leq2\left(
1-1/\omega\right)  \left(  1-1/\omega\right)  n^{2}/2=\left(  \left(
1-1/\omega\right)  n\right)  ^{2},
\]
and so (\ref{bnd6}) implies (\ref{bnd4}).

As shown in \cite{Nik02}, inequality (\ref{bnd6}) follows from the celebrated
result of Motzkin and Straus \cite{MoSt65}:\medskip

\emph{Let }$G$\emph{ be a graph of order }$n$\emph{ with cliques number
}$\omega\left(  G\right)  =\omega.$\emph{ If }$\left(  x_{1},\ldots
,x_{n}\right)  $\emph{ is a vector with nonnegative entries, then}%
\begin{equation}
\sum_{uv\in E\left(  G\right)  }x_{u}x_{v}\leq\frac{\omega-1}{2\omega}\left(
\sum_{u\in V\left(  G\right)  }x_{u}\right)  ^{2}. \label{inMM}%
\end{equation}

On the other hand, this result follows in turn from the concise Tur\'{a}n
theorem, as shown in \cite{Nik06a}. The implications
\[
\text{(\ref{bnd6})}\Longrightarrow\text{(\ref{conTT})}\Longrightarrow
\text{MS}\Longrightarrow\text{(\ref{bnd6})}%
\]
justify regarding inequality (\ref{bnd6}) as a flexible spectral form of the
concise Tur\'{a}n theorem.

Next we extend Theorem \ref{Nik02} in a somewhat unexpected direction. Recall
that, a $k$\emph{-walk} in a graph $G$ is a sequence of vertices
$v_{1},...,v_{k}$ of $G$ such that $v_{i}$ is adjacent to $v_{i+1}$ for
$i=1,...,k-1;$ write $w_{k}\left(  G\right)  $ for the number of $k$-walks in
$G.$ Observing that $2e\left(  G\right)  =w_{2}\left(  G\right)  ,$ we see
that the following theorem, given in \cite{Nik06} , extends inequality
(\ref{bnd6}).

\begin{theorem}
\label{thUB}If $r\geq1$ and $G$ is a graph with clique number $\omega\left(
G\right)  =\omega,$ then%
\begin{equation}
\mu^{r}\left(  G\right)  \leq\left(  1-1/\omega\right)  w_{r}\left(  G\right)
. \label{upb}%
\end{equation}
Suppose that $G$ has no isolated vertices and equality holds for some $r\geq1$.

(i) If $r=1$, then $G$ is a regular complete $\omega$-partite graph.

(ii) If $r\geq2$ and $\omega>2$, then $G$ is a regular complete $\omega
$-partite graph.

(iii) If $r\geq2$ and $\omega=2$, then $G$ is a complete bipartite graph, and
if $r$ is odd, then $G$ is regular.
\end{theorem}

It is somewhat surprising that for $r\geq2$ the number of vertices of $G$ is
not relevant in this theorem.

\subsubsection*{A precise spectral Tur\'{a}n theorem}

Since Wilf's inequality (\ref{bnd4}) becomes equality only when $\omega$
divides $n,$ one can expect that some fine tuning is still possible. Indeed,
in \cite{Nik07} we sharpened inequality (\ref{bnd4}), bringing it the closest
possible to the Tur\'{a}n theorem.

\begin{theorem}
If $G$ is a graph of order $n$ with no complete subgraph of order $r+1$, then
$\mu\left(  G\right)  \leq\mu\left(  T_{r}\left(  n\right)  \right)  .$
Equality holds if and only if $G=T_{r}\left(  n\right)  .$
\end{theorem}

Here is an equivalent, shorter form of this statement:\emph{ If }%
$G\in\mathcal{G}\left(  n\right)  $\emph{ and }$\omega\left(  G\right)
=\omega,$ \emph{then }$\mu\left(  G\right)  <\mu\left(  T_{\omega}\left(
n\right)  \right)  $\emph{ unless }$G=T_{\omega}\left(  n\right)  .$

Note also that $\mu\left(  T_{2}\left(  n\right)  \right)  =\sqrt{\left\lfloor
n^{2}/4\right\rfloor };$ for $\omega\geq3$ there is also a closed expression
for $\mu\left(  T_{\omega}\left(  n\right)  \right)  ,$ but it is somewhat cumbersome.

\subsubsection*{Spectral radius and independence number}

One wonders if there is a theorem about the independence number $\alpha\left(
G\right)  $, similar to the Tur\'{a}n theorem. One obvious answer is obtained
by restating the concise Tur\'{a}n theorem in complementary terms%
\[
2e\left(  G\right)  \geq n^{2}/\alpha\left(  G\right)  -n,
\]
which immediately implies that $\mu\left(  G\right)  \geq n/\alpha\left(
G\right)  -1$ as well. Note that here the spectral statement follows from the
conventional one. However, a proof by induction on $\alpha$ gives the
following sharper result.

\begin{theorem}
\label{in}If $G\in G\left(  n\right)  $ and $\alpha\left(  G\right)  =\alpha$,
then $\mu\left(  G\right)  \geq\left\lceil n/\alpha\right\rceil -1.$
\end{theorem}

In a different direction, for connected graphs and some special values of
$\alpha,$ more specific results have been proved in \cite{XHSZ09}.

Also, by the well-known inequality $q_{1}\left(  G\right)  \geq2\mu_{1}\left(
G\right)  ,$ Theorem \ref{in} proves Conjecture 27 from \cite{HaLu10}.

\subsection{A spectral Erd\H{o}s-Stone-Bollob\'{a}s theorem}

Having seen various spectral forms of the Tur\'{a}n theorem, one can expect
that many other results that surround it can be cast in spectral form as well;
and this is indeed the case. The following theorem, given in \cite{Nik09f}, is
the spectral analog of the Erd\H{o}s-Stone-Bollob\'{a}s theorem, more
precisely, of Theorem \ref{ESnew}.

\begin{theorem}
\label{genZ}Let $r\geq3,$ let $c$ and $n$ be such that
\[
0<c<1/r!,\text{ \ \ }n\geq\exp\left(  \left(  r^{r}/c\right)  ^{r}\right)  ,
\]
and let $G$ be a graph of order $n$. If
\begin{equation}
\mu\left(  G\right)  \geq\left(  1-1/\left(  r-1\right)  +c\right)  n,
\label{bnd}%
\end{equation}
then $G$ contains a $K_{r}\left(  s;t\right)  $ with $s\geq\left\lfloor
\left(  c/r^{r}\right)  ^{r}\log n\right\rfloor $ and$\ t>n^{1-c^{r-1}}.$
\end{theorem}

Let us emphasize that the functionality of Theorem \ref{ESnew} is entirely
preserved: in particular, $c$ may depend on $n,$ e.g., letting $c=1/\log\log
n,$ the conclusion is meaningful for sufficiently large $n.$

Since $\mu\left(  G\right)  \geq2e\left(  G\right)  /v\left(  G\right)  ,$
Theorem \ref{genZ}, in fact, implies Theorem \ref{ESnew}. Other lower bounds
on $\mu\left(  G\right)  ,$ such as those given in \cite{Nik06}, imply other
new versions of this theorem.

Suppose that $c$ is a sufficiently small constant. Choosing randomly a graph
$G$ of order $n$ with $\left\lceil \left(  1-1/\left(  r-1\right)  +2c\right)
n^{2}/2\right\rceil $ edges, we have $\mu\left(  G\right)  \geq\left(
1-1/\left(  r-1\right)  +c\right)  n$, but $G$ contains no $K_{2}\left(
\left\lfloor C\log\right\rfloor ,\left\lfloor C\log n\right\rfloor \right)  $
for some $C>0,$ independent of $n.$ Hence, for constant $c,$ Theorem
\ref{genZ} is best possible up to a constant factor.

We close this topic with a consequence of Theorem \ref{genZ}, given in
\cite{Nik09f}, that solves asymptotically the following general extremal
problem: \emph{Given a family $\mathcal{F}$ of nonbipartite forbidden
subgraphs}$,$\emph{ what is the maximum spectral radius of a graph of order
}$n$\emph{ containing no member of $\mathcal{F}$}.

\begin{theorem}
\label{co1}Let $r\geq3$ and let $F_{1},F_{2},\ldots$ be $r$-partite graphs
satisfying $v\left(  F_{n}\right)  =o\left(  \log n\right)  .$ Then%
\begin{equation}
\max\left\{  \mu\left(  G\right)  :G\in\mathcal{G}\left(  n\right)  \text{ and
}F_{n}\nsubseteq G\right\}  =\left(  1-\frac{1}{r-1}\right)  n+o\left(
n\right)  . \label{b1}%
\end{equation}

\end{theorem}

It is likely that in the setup of Theorem \ref{co1} the condition $v\left(
F_{n}\right)  =o\left(  \log n\right)  $ can be sharpened.

\subsection{Saturation problems}

The precise spectral Tur\'{a}n theorem implies that if $G$ is a graph of order
$n$ with $\mu\left(  G\right)  >\mu\left(  T_{r}\left(  n\right)  \right)  ,$
then $G$ contains a $K_{r+1}.$ Since this setting is analogous to the case
when $e\left(  G\right)  >t_{r}\left(  n\right)  ,$ one would expect much
larger supergraphs of $K_{r+1}$. In fact, as we shall see, all results from
Subsection \ref{sat} have their spectral analogs. In addition, it is not
difficult to show that if $G$ is a graph of order $n,$ then the inequality
$e\left(  G\right)  >e\left(  T_{r}\left(  n\right)  \right)  $ implies the
inequality $\mu\left(  G\right)  >\mu\left(  T_{r}\left(  n\right)  \right)
.$ Therefore, the spectral theorems below imply the corresponding nonspectral
extremal results, albeit with somewhat narrower ranges of the parameters.

We start with the spectral analog of Theorem \ref{Th1}, given in \cite{Nik09b}.

\begin{theorem}
\label{spN2} Let $r\geq2,$ let $c$ and $n$ be such that%
\[
0<c\leq r^{-\left(  2r+9\right)  \left(  r+1\right)  },\text{ \ \ }n\geq
\exp\left(  2/c\right)  ,
\]
and let $G$ be a graph of order $n.$ If $\mu\left(  G\right)  >\mu\left(
T_{r}\left(  n\right)  \right)  ,$ then $G$ contains a
\[
K_{r}^{+}\left(  \left\lfloor c\log n\right\rfloor ;\left\lceil n^{1-\sqrt{c}%
}\right\rceil \right)  .
\]

\end{theorem}

Theorem \ref{spN2} is essentially best possible since for every $\varepsilon
>0,$ choosing randomly a graph $G$ of order $n$ with $e\left(  G\right)
=\left\lceil \left(  1-\varepsilon\right)  n^{2}/2\right\rceil ,$ we see that
$\mu\left(  G\right)  >\left(  1-\varepsilon\right)  n,$ but $G$ contains no
$K_{2}\left(  \left\lfloor c\log n\right\rfloor \right)  $ for some $c>0,$
independent of $n.$

The theorem corresponding to Theorem \ref{tJ} is given in \cite{Nik09b}. We
state it here in a somewhat refined form.

\begin{theorem}
\label{spN1}Let $r\geq2,$ $n>r^{15},$ and let $G$ be a graph of order $n.$ If
$\mu\left(  G\right)  \geq\mu\left(  T_{r}\left(  n\right)  \right)  ,$ then
\[
js_{r+1}\left(  G\right)  >n^{r-1}/r^{2r+4}%
\]
unless $G=T_{r}\left(  n\right)  .$
\end{theorem}

Theorem \ref{spN1} and its stability complement Theorem \ref{spN1.2} are
crucial in the proof of several other spectral extremal results.

It is easy to see the Tur\'{a}n graph $T_{2}\left(  n\right)  $ contains no
odd cycles and that $\mu\left(  T_{2}\left(  n\right)  \right)  =\sqrt
{\left\lfloor n^{2}/4\right\rfloor }.$ Hence the following theorem gives a
sharp spectral condition for the existence of odd cycles.

\begin{theorem}
\label{spN3}Let $G$ be a graph of sufficiently large order $n.$ If $\mu\left(
G\right)  >\sqrt{\left\lfloor n^{2}/4\right\rfloor },$ then $G$ contains a
cycle of length $t$ for every $t\leq n/320.$
\end{theorem}

This theorem, given in \cite{Nik09g}, is motivated by the following result of
Bollob\'{a}s (\cite{Bol78}, p. 150): \emph{if }$G$\emph{ is a graph of order
}$n$\emph{ with }$e\left(  G\right)  >\left\lfloor n^{2}/4\right\rfloor
$\emph{, then }$G$\emph{ contains a cycle of length }$t$\emph{ for every
}$t=3,\ldots,\left\lceil n/2\right\rceil .$

\subsection{Stability problems}

We shall show that most stability results from Subsection \ref{stab} have
their spectral analogs. However, we could not find spectral analogs of the
stability problems that involve minimum degree (Subsection \ref{3free}).

We first state a general spectral stability result, and then two stronger
versions for specific graphs. We give here only Theorems \ref{tstab} and
\ref{spN1.2} since they are important for various applications, but our
machinery helped to deduce many others of somewhat lesser importance, and they
can be found in \cite{Nik09b} and \cite{Nik09g}.

The following analog of Theorem \ref{tStab} was given in \cite{Nik09b}.

\begin{theorem}
\label{stStab}Let $r\geq2,$ let $c,$ $\varepsilon$ and $n$ be such that
\[
0<c<r^{-8\left(  r+21\right)  \left(  r+1\right)  },\text{\ \ \ }%
0<\varepsilon<2^{-36}r^{-24},\text{ \ \ }n>\exp\left(  1/c\right)  ,
\]
and let $G$ be a graph of order $n.$ If $\mu\left(  G\right)  >\left(
1-1/r-\varepsilon\right)  n$, then one of the following statements holds:

(a) $G$ contains a $K_{r+1}\left(  \left\lfloor c\log\right\rfloor
;\left\lceil n^{1-\sqrt{c}}\right\rceil \right)  ;$

(b) $G$ differs from $T_{r}\left(  n\right)  $ in fewer than $\left(
\varepsilon^{1/4}+c^{1/\left(  8r+8\right)  }\right)  n^{2}$ edges.
\end{theorem}

The proofs of Theorem \ref{tStab} and \ref{stStab}, given in \cite{Nik09b}
illustrate the isomorphism between the sets of tools developed for the
spectral and nonspectral problems. The texts of the two proofs are almost
identical, while the differences come from the use of different tools. We
refer the reader to Section \ref{spf} for more details.

The next two theorems are crucial for several applications. The first one,
proved in \cite{BoNi07}, is a spectral equivalent of Theorem \ref{stabK}.

\begin{theorem}
\label{tstab}Let $r\geq2$ and $0\leq\varepsilon\leq2^{-10}r^{-6},$ and let $G$
be a $K_{r+1}$-free graph of order $n.$ If
\begin{equation}
\mu\left(  G\right)  \geq\left(  1-1/r-\varepsilon\right)  n, \label{reqmu}%
\end{equation}
then $G$ contains an induced $r$-partite graph $H$ of order at least $\left(
1-3\alpha^{1/3}\right)  n$ and minimum degree%
\[
\delta\left(  H\right)  >\left(  1-1/r-6\varepsilon^{1/3}\right)  n.
\]

\end{theorem}

Finally, we have a spectral stability theorem for large joints, proved in
\cite{Nik09b}.

\begin{theorem}
\label{spN1.2}Let $r\geq2,$ let $\varepsilon$ and $n$ be such that
\[
0<\varepsilon<2^{-10}r^{-6},\text{ \ \ }n\geq r^{20},
\]
and let $G$ be a graph of order $n.$ If $\mu\left(  G\right)  >\left(
1-1/r-\varepsilon\right)  n,$ then $G$ satisfies one of the conditions:

(a) $js_{r+1}\left(  G\right)  >n^{r-1}/r^{2r+5};$

(b) $G$ contains an induced $r$-partite subgraph $H$ of order at least
$\left(  1-4\varepsilon^{1/3}\right)  n$ with minimum degree $\delta\left(
H\right)  >\left(  1-1/r-7\varepsilon^{1/3}\right)  n.$
\end{theorem}

\subsection{The Zarankiewicz problem}

What is the maximum spectral radius of a graph of order $n$ with no $K_{s,t}?$
This is a spectral version of the famous Zarankiewicz problem: \emph{what is
the maximum number of edges in a graph of order }$n$\emph{ with no }$K_{s,t}?$
Except for few cases, no complete solution to either of these problems is
known. For instance, Babai and Guiduli \cite{BaGu09} have shown that
\[
\mu\leq\left(  \left(  s-1\right)  ^{1/t}+o\left(  1\right)  \right)
n^{1-1/t}.
\]
Using a different method, in \cite{Nik10b} we improved this result as follows:

\begin{theorem}
\label{tZar}Let $s\geq t\geq2,$ and let $G$ be a $K_{s,t}$-free graph of order
$n.$

(i) If $t=2,$ then%
\begin{equation}
\mu\left(  G\right)  \leq1/2+\sqrt{\left(  s-1\right)  \left(  n-1\right)
+1/4}. \label{in0}%
\end{equation}

(ii) If $t\geq3,$ then%
\begin{equation}
\mu\left(  G\right)  \leq\left(  s-t+1\right)  ^{1/t}n^{1-1/t}+\left(
t-1\right)  n^{1-2/t}+t-2. \label{in1}%
\end{equation}

\end{theorem}

On the other hand, in view of the inequality $2e\left(  G\right)  \leq
\mu\left(  G\right)  n,$ we see that if $G$ is a $K_{s,t}$-free graph of order
$n,$ then
\begin{equation}
e\left(  G\right)  \leq\frac{1}{2}\left(  s-t+1\right)  ^{1/t}n^{2-1/t}%
+\frac{1}{2}\left(  t-1\right)  n^{2-2/t}+\frac{1}{2}\left(  t-2\right)  n.
\label{in2}%
\end{equation}
This is a slight improvement of a result of F\"{u}redi \cite{Fur96} and this
seems the best known bound on $e\left(  G\right)  $ so far.

For some values of $s$ and $t$ the bounds given by (\ref{in0}) and (\ref{in1})
are tight as we now demonstrate.

\subsubsection*{\textbf{The case }$t=2$}

For $s=t=2$ inequality (\ref{in0}) shows that every $K_{2,2}$-free graph $G$
of order $n$ satisfies
\[
\mu\left(  G\right)  \leq1/2+\sqrt{n-3/4}.
\]
This bound is tight because equality holds for the friendship graph (a
collection of triangles sharing a single common vertex).

Also, Erd\H{o}s-Renyi \cite{ErRe62} showed that if $q$ is a prime power, the
polarity graph $ER_{q}$ is a $K_{2,2}$-free graph of order $n=q^{2}+q+1$ and
$q\left(  q+1\right)  ^{2}/2$ edges. Thus, its spectral radius $\mu\left(
ER_{q}\right)  $ satisfies
\[
\mu\left(  ER_{q}\right)  \geq\frac{q^{3}+2q^{2}+q}{q^{2}+q+1}>q+1-\frac{1}%
{q}=1/2+\sqrt{n-3/4}-\frac{1}{\sqrt{n}-1},
\]
which is also close to the upper bound.

For $s>2,$ equality in (\ref{in0}) is attained when $G$ is a strongly regular
graph in which every two vertices have exactly $s-1$ common neighbors. There
are examples of strongly regular graphs of this type; here is a small
selection from Gordon Royle's webpage:%

\[%
\begin{tabular}
[c]{||c|c|c||}\hline\hline
\emph{s} & \emph{n} & $\mu\left(  G\right)  $\\
$3$ & $45$ & $12$\\
$4$ & $96$ & $20$\\
$5$ & $175$ & $30$\\
$6$ & $36$ & $15$\\\hline\hline
\end{tabular}
\ \
\]
We are not aware whether there are infinitely many strongly regular graphs in
which every two vertices have the same number of common neighbors. However,
F\"{u}redi \cite{Fur96a} has shown that for any $n$ there exists a $K_{s,2}%
$-free graph $G_{n}$ of order $n$ such that%
\[
e\left(  G_{n}\right)  \geq\frac{1}{2}n\sqrt{sn}+O\left(  n^{4/3}\right)  ,
\]
and so,%
\[
\mu\left(  G_{n}\right)  \geq\sqrt{sn}+O\left(  n^{1/3}\right)  ;
\]
thus (\ref{in0}) is tight up to low order terms.

\subsubsection*{\textbf{The case} $s=t=3$}

The bound (\ref{in1}) implies that if $G$ is a $K_{3,3}$-free graph of order
$n,$ then
\[
\mu\left(  G\right)  \leq n^{2/3}+2n^{1/3}+1\text{\textbf{.}}%
\]

On the other hand, a construction due to Alon, R\`{o}nyai and Szab\`{o}
\cite{ARS99} implies that for all $n=q^{3}-q^{2},$ where $q$ is a prime power,
there exists a $K_{3,3}$-free graph $G_{n}$ of order $n$ with
\[
\mu\left(  G_{n}\right)  \geq n^{2/3}+\frac{2}{3}n^{1/3}+C
\]
for some constant $C>0.$ Thus, the bound (\ref{in1}) is asymptotically tight
for $s=t=3.$ The same conclusion can be obtained from Brown's construction of
$K_{3,3}$-free graphs \cite{Bro66}.

\subsubsection*{\textbf{The general case}}

As proved in \cite{ARS99}, there exists $c>0$ such that for all $t\geq2$ and
$s\geq\left(  t-1\right)  !+1$, there is a $K_{s,t}$-free graph $G_{n}$ of
order $n$ with
\[
e\left(  G_{n}\right)  \geq\frac{1}{2}n^{2-1/t}+O\left(  n^{2-1/t-c}\right)
.
\]
Hence, for such $s$ and $t$ we have%
\[
\mu\left(  G\right)  \geq n^{1-1/t}+O\left(  n^{1-1/t-c}\right)  ;
\]
thus, the bound (\ref{in1}) and also the bound of Babai and Guiduli give the
correct order of the main term.

\subsection{Paths and cycles}

We give now some results about the maximum spectral radius of graphs of order
$n$ without paths or cycles of specified length. Writing $C_{k}$ and $P_{k}$
for the cycle and path of order $k$, let us define the functions%
\begin{align*}
f_{l}\left(  n\right)   &  =\max\left\{  \mu\left(  G\right)  :G\in
\mathcal{G}\left(  n\right)  \text{ and }C_{l}\nsubseteq G\right\}  ;\\
g_{l}\left(  n\right)   &  =\max\left\{  \mu\left(  G\right)  :G\in
\mathcal{G}\left(  n\right)  \text{ and }C_{l}\nsubseteq G,\text{ and }%
C_{l+1}\nsubseteq G\right\}  ;\\
h_{l}\left(  n\right)   &  =\max\left\{  \mu\left(  G\right)  :G\in
\mathcal{G}\left(  n\right)  \text{ and }P_{l}\nsubseteq G\right\}  .
\end{align*}

For these functions we shall show below some exact expressions or at least
good asymptotics. It should be noted that except for $f_{l}\left(  n\right)  $
when $l$ is odd, these questions are quite different from their nonspectral analogs.

The lower bounds on $f_{2l}\left(  n\right)  ,$ $g_{l}\left(  n\right)  $ and
$h_{l}\left(  n\right)  $ are given by two families of graphs, which for
sufficiently large $n$ give the exact values of $h_{l}\left(  n\right)  $, and
perhaps also of $f_{2l}\left(  n\right)  $ and $g_{l}\left(  n\right)
.$\medskip

\emph{Suppose that }$1\leq k<n$\emph{.}

\emph{(1) Let }$S_{n,k}$\emph{ be the graph obtained by joining every vertex
of a complete graph of order }$k$\emph{ to every vertex of an independent set
of order }$n-k,$\emph{ that is, }$S_{n,k}=K_{k}\vee\overline{K}_{n-k};$

\emph{(2) Let }$S_{n,k}^{+}$\emph{ be the graph obtained by adding one edge
within the independent set of }$S_{n,k}.$\medskip

Note that $P_{l+1}\nsubseteq S_{n,k}$ and $C_{l}\nsubseteq S_{n,k}$ for
$l\geq2k+1.$ Likewise, $P_{l+1}\nsubseteq S_{n,k}$ and $C_{l}\nsubseteq
S_{n,k}$ for $l\geq2k+2.$

Therefore,
\begin{align*}
h_{2k}\left(  n\right)   &  \geq\mu\left(  S_{n,k}\right)  =\left(
k-1\right)  /2+\sqrt{kn-\left(  3k^{2}+2k-1\right)  /4},\\
h_{2k+1}\left(  n\right)   &  \geq\mu\left(  S_{n,k}^{+}\right)  =\left(
k-1\right)  /2+\sqrt{kn-\left(  3k^{2}+2k-1\right)  /4}+1/n+O\left(
n^{-3/2}\right)  ,\\
g_{2k}\left(  n\right)   &  \geq\mu\left(  S_{n,k}\right)  =\left(
k-1\right)  /2+\sqrt{kn-\left(  3k^{2}+2k-1\right)  /4},\\
g_{2k+1}\left(  n\right)   &  \geq\mu\left(  S_{n,k}^{+}\right)  =\left(
k-1\right)  /2+\sqrt{kn-\left(  3k^{2}+2k-1\right)  /4}+1/n+O\left(
n^{-3/2}\right)  ,\\
f_{2k+2}\left(  n\right)   &  \geq\mu\left(  S_{n,k}^{+}\right)  =\left(
k-1\right)  /2+\sqrt{kn-\left(  3k^{2}+2k-1\right)  /4}+1/n+O\left(
n^{-3/2}\right)  .
\end{align*}
Below we shall give also rather close upper bounds for these functions.

\subsubsection*{Forbidden odd cycle}

In view of Theorem \ref{spN3}, we find that if $l$ is odd and $n>320l,$ then%
\[
f_{l}\left(  n\right)  =\sqrt{\left\lfloor n^{2}/4\right\rfloor }.
\]
The smallest ratio $n/l$ for which this equation is still valid is not known.

Clearly, for odd $l$ we have $f_{l}\left(  n\right)  \sim n/2,$ which is in
sharp contrast to the value of $f_{l}\left(  n\right)  $ for even $l.$

\subsubsection*{Forbidden cycle $C_{4}$}

The value of $f_{4}\left(  n\right)  $ is essentially determined in
\cite{Nik07}:\medskip

\emph{Let }$G$\emph{ be a graph of order }$n$\emph{ with }$\mu\left(
G\right)  =\mu$\emph{. If }$C_{4}\nsubseteq G,$\emph{ then}%
\begin{equation}
\mu^{2}-\mu\leq n-1. \label{c4in}%
\end{equation}
\emph{Equality holds if and only if every two vertices of }$G$\emph{ have
exactly one common neighbor, i.e., when }$G$\emph{ is the friendship
graph.}\medskip

An easy calculation implies that\emph{ }%
\[
f_{4}\left(  n\right)  =1/2+\sqrt{n-3/4}+O\left(  1/n\right)  ,
\]
where for odd $n$ the $O\left(  1/n\right)  $ term is zero. Finding the
precise value of $f_{4}\left(  n\right)  $ for even $n$ is an open problem.

Here is a considerably more involved bound on the spectral radius of a $C_{4}%
$-free graph of given size, given in \cite{Nik09e}.

\begin{theorem}
\label{tC1}Let $m\geq9$ and $G$ be a graph with $m$ edges. If $\mu\left(
G\right)  >\sqrt{m},$ then $G$ has a $4$-cycle.
\end{theorem}

This theorem is tight, for all stars are $C_{4}$-free graphs with $\mu\left(
G\right)  =\sqrt{m}.$ Also, let $S_{n,1}$ be the star of order $n$ with an
edge within its independent set. The graph $S_{n,1}$ is $C_{4}$-free and has
$n$ edges, but $\mu\left(  G\right)  >\sqrt{n}$ for $4\leq n\leq8,$ while
$\mu\left(  S_{9,1}\right)  =3$.

\subsubsection*{Forbidden cycle $C_{2k}$}

The inequality (\ref{c4in}) can be generalized for arbitrary even cycles as
follows: \emph{if }$C_{2k+2}\nsubseteq G,$\emph{ then}%
\[
\mu^{2}-\left(  k-1\right)  \mu\leq k\left(  n-1\right)  .
\]
In fact, a slightly stronger assertion was proved in \cite{Nik10a}.

\begin{theorem}
\label{tP2} Let $k\geq1$ and $G$ be a graph of order $n.$ If
\[
\mu\left(  G\right)  >k/2+\sqrt{kn+\left(  k^{2}-4k\right)  /4},
\]
then $C_{2l+2}\subset G$ for every $l=1,\ldots,k$.
\end{theorem}

In view of the graph $S_{n,k}^{+},$ Theorem \ref{tP2} implies that%
\begin{equation}
\left(  k-1\right)  /2+\sqrt{kn}+o\left(  n\right)  \leq f_{2k+2}\left(
n\right)  \leq k/2+\sqrt{kn}+o\left(  n\right)  . \label{cyc}%
\end{equation}
The exact value of $f_{2k+2}\left(  n\right)  $ is not known for $k\geq2,$ and
finding this value seems a challenge. Nevertheless, the precision of
(\ref{cyc}) is somewhat surprising, given that the asymptotics of the maximum
number of edges in $C_{2k+2}$-free graphs of order $n$ is not known for
$k\geq2$.

\subsubsection*{Forbidden pair of cycles $\left\{  C_{2k},C_{2k+1}\right\}  $}

Let us consider now the function $g_{l}\left(  n\right)  $.\ To begin with,
Favaron, Mah\'{e}o, and Sacl\'{e} \cite{FMS93} showed that if a graph $G$ of
order $n$ contains neither $C_{3}$ nor $C_{4},$ then $\mu\left(  G\right)
\leq\sqrt{n-1}.$ Since the star of order $n$ has no cycles and its spectral
radius is $\sqrt{n-1},$ we see that
\[
g_{3}\left(  n\right)  =\sqrt{n-1}\emph{.}%
\]
We do not know the exact value of $g_{l}\left(  n\right)  $ for $l>3,$ but we
have the following theorem from \cite{Nik10a}.

\begin{theorem}
\label{tP3}Let $k\geq1$ and $G$ be a graph of order $n.$ If%
\[
\mu\left(  G\right)  >\left(  k-1\right)  /2+\sqrt{kn+\left(  k+1\right)
^{2}/4},
\]
then $C_{2k+1}\subset G$ or $C_{2k+2}\subset G.$
\end{theorem}

Theorem \ref{tP3}, together with the graphs $S_{n,k}$ and $S_{n,k}^{+}$,
gives
\begin{align*}
\left(  k-1\right)  /2+\sqrt{kn}+o\left(  n\right)   &  \leq g_{2k+1}\left(
n\right)  \leq k/2+\sqrt{kn}+o\left(  n\right)  ,\\
g_{2k}\left(  n\right)   &  =\left(  k-1\right)  /2+\sqrt{kn}+\Theta\left(
n^{-1/2}\right)  .
\end{align*}

\subsubsection*{Forbidden path $P_{k}$}

The function $h_{k}\left(  n\right)  $ is completely known for large $n.$ As
proved in \cite{Nik10a}:

\begin{theorem}
\label{tP1} Let $k\geq1,$ $n\geq2^{4k}$ and let $G$ be a graph of order $n.$

(i) If $\mu\left(  G\right)  \geq\mu\left(  S_{n,k}\right)  ,$ then $G$
contains a $P_{2k+2}$ unless $G=S_{n,k}.$

(ii) If $\mu\left(  G\right)  \geq\mu\left(  S_{n,k}^{+}\right)  ,$ then $G$
contains a $P_{2k+3}$ unless $G=S_{n,k}^{+}.$
\end{theorem}

Theorem \ref{tP1}, together with the graphs $S_{n,k}$ and $S_{n,k}^{+}$,
implies that for every $k\geq1$ and $n\geq2^{4k},$ we have
\begin{align*}
h_{2k}\left(  n\right)   &  =\mu\left(  S_{n,k}\right)  =\left(  k-1\right)
/2+\sqrt{kn-\left(  3k^{2}+2k-1\right)  /4},\\
h_{2k+1}\left(  n\right)   &  =\mu\left(  S_{n,k}^{+}\right) \\
&  =\left(  k-1\right)  /2+\sqrt{kn-\left(  3k^{2}+2k-1\right)  /4}%
+1/n+O\left(  n^{-3/2}\right)  .
\end{align*}

\subsection{Hamilton paths and cycles}

In \cite{Ore60}, Ore found the following sufficient condition for the
existence of Hamilton paths and cycles.

\begin{theorem}
\label{f1}Let $G$ be a graph of order $n$. If
\begin{equation}
e\left(  G\right)  \geq\binom{n-1}{2} \label{siz}%
\end{equation}
then $G$ contains a Hamiltonian path unless $G=K_{n-1}+K_{1}.$ If the
inequality (\ref{siz}) is strict, then $G$ contains a Hamiltonian cycle unless
$G=K_{n-1}+e.$
\end{theorem}

In the line above and further, $K_{n-1}+e$ denotes the complete graph
$K_{n-1}$ with a pendent edge.

Recently, Fiedler and Nikiforov \cite{FiNi10} deduced a spectral version of
Ore's result.

\begin{theorem}
\label{tHam1}Let $G$ be a graph of order $n.$ If
\begin{equation}
\mu\left(  G\right)  \geq n-2, \label{mu}%
\end{equation}
then $G$ contains a Hamiltonian path unless $G=K_{n-1}+K_{1}.$ If the
inequality (\ref{mu}) is strict, then $G$ contains a Hamiltonian cycle unless
$G=K_{n-1}+e.$
\end{theorem}

A subtler spectral condition for Hamiltonicity was obtained using the spectral
radius of the complement of a graph.

\begin{theorem}
\label{tHam2}Let $G$ be a graph of order $n$ and $\mu\left(  \overline
{G}\right)  $ be the spectral radius of its complement. If
\[
\mu\left(  \overline{G}\right)  \leq\sqrt{n-1},
\]
then $G$ contains a Hamiltonian path unless $G=K_{n-1}+K_{1}.$ If
\[
\mu\left(  \overline{G}\right)  \leq\sqrt{n-2},
\]
then $G$ contains a Hamiltonian cycle unless $G=K_{n-1}+e.$
\end{theorem}

Zhou \cite{Zho10}, adopting the same technique, proved a similar result for
the signless Laplacian, which has been subsequently refined in \cite{Nik10d}
to the following one:

\begin{theorem}
\label{tQ}Let $G$ be a graph of order $n$ and $q\left(  \overline{G}\right)  $
be the spectral radius of the $Q$-matrix of its complement.

(i) If
\begin{equation}
q\left(  \overline{G}\right)  \leq n, \label{cz1}%
\end{equation}
then $G$ contains a Hamiltonian path unless $G$ is the union of two disjoint
complete graphs or $n$ is even and $G=K_{n/2-1,n/2+1}.$

(ii) If
\begin{equation}
q\left(  \overline{G}\right)  \leq n-1, \label{cz2}%
\end{equation}
then $G$ contains a Hamiltonian cycle unless $G$ is a union of two complete
graphs with a single common vertex or $n$ is odd and $G=K_{\left\lfloor
n/2\right\rfloor ,\left\lceil n/2\right\rceil }.$
\end{theorem}

Note that if the inequality in (\ref{cz1}) or (\ref{cz2}) is strict, then the
corresponding conclusion holds with no exception. Also, as it turns out,
Theorem \ref{tQ} considerably strengthens the classical degree conditions for
Hamiltonicity by Ore \cite{Ore60}.

\subsection{Clique number and eigenvalues}

If a triangle-free graph is sufficiently dense, then it contains large
independent sets and the modulus of its smallest eigenvalue cannot be very
small. A more general statement of such type has been proved in \cite{BoNi04a}
for graphs of bounded clique number. Somewhat later, the following explicit
dependence was found in \cite{Nik06}.

\begin{theorem}
\label{minK}If $G\in\mathcal{G}\left(  n,m\right)  $ and $\omega\left(
G\right)  =\omega$, then%
\begin{equation}
\mu_{n}\left(  G\right)  <-\frac{2}{\omega}\left(  \frac{2m}{n^{2}}\right)
^{\omega}n. \label{bnd1}%
\end{equation}

\end{theorem}

Inequality (\ref{bnd1}) captures pretty well the situation in dense graphs,
that is, if $G$ is a dense graph with $\mu_{n}\left(  G\right)  =O\left(
n^{1-c}\right)  $ for some $c\in\left(  0,1/2\right)  ,$ then $G$ contains
cliques of order $\Omega\left(  \log n\right)  $.

Moreover, as shown in \cite{Nik06}, inequality (\ref{bnd1}) is tight up to a
constant factor for several classes of sparse graphs, but complete
investigation of this issue seems difficult.

In \cite{Nik09d}, inequality (\ref{bnd1}) was used to derive a lower bound on
$\alpha\left(  G\right)  ,$ thus giving other cases of tightness.

\begin{theorem}
\label{thm1} Let $G\in\mathcal{G}\left(  n,m\right)  ,$ $d=2m/n,\ $and
$\tau=\left\vert \mu_{n}\left(  \overline{G}\right)  \right\vert .$ If
$d\geq2,$ then
\[
\alpha\left(  G\right)  >\left(  \frac{n}{d+1}-1\right)  \left(  \log
\frac{d+1}{\tau}-\log\log\left(  d+1\right)  \right)  .
\]

\end{theorem}

Inequality (\ref{bnd1}) is concise, but it is difficult to use because the
right-hand side is exponential in $\omega\left(  G\right)  $. The following
two somewhat simpler bounds, given in \cite{Nik07b}, stem from Tur\'{a}n's
theorem and some inequalities that will be given in the next subsection.

\begin{theorem}
\label{thm3} Let $G\in\mathcal{G}\left(  n,m\right)  ,$ $d=2m/n,$ and
$\omega\left(  G\right)  =\omega.$ Then
\[
\omega\geq1+\frac{dn}{\left(  n-d\right)  \left(  d-\mu_{n}\left(  G\right)
\right)  }.
\]
Equality holds if and only if $G$ is a complete regular $\omega$-partite graph.
\end{theorem}

Similar inequalities \cite{Nik07} exist also for the Laplacian eigenvalues.

\begin{theorem}
Let $G=\mathcal{G}\left(  n,m\right)  ,$ $d=2m/n$ and $\omega\left(  G\right)
=\omega.$ Then
\[
\omega\geq1+\frac{dn}{\lambda_{n}\left(  G\right)  \left(  n-d\right)  },
\]
with equality holding if and only if $G$ is a regular complete $\omega
$-partite graph.

Also,
\[
\alpha\left(  G\right)  \geq1+\frac{\left(  n-1-d\right)  n}{\left(
n-\lambda_{2}\left(  G\right)  \right)  \left(  1+d\right)  },
\]
with equality holding if and only if $G$ is the union of $\alpha\left(
G\right)  $ disjoint cliques of equal order.
\end{theorem}

Note that both bounds in the last theorem imply the concise Tur\'{a}n theorem.

\subsection{Number of cliques and eigenvalues}

It turns out that the numbers of various cliques of a graph are closely
related to its most important eigenvalues. Bollob\'{a}s and Nikiforov
\cite{BoNi07} proved the following chain of inequalities, which were useful on
several occasions.

\begin{theorem}
\label{le3mu}Let $G$ be a graph with $\omega\left(  G\right)  =\omega\geq2$
and $\mu\left(  G\right)  =\mu.$ For every $r=2,\ldots,\omega,$%
\[
\mu^{r+1}\leq\left(  r+1\right)  k_{r+1}\left(  G\right)  +\sum_{s=2}%
^{r}\left(  s-1\right)  k_{s}\left(  G\right)  \mu^{r+1-s}.
\]

\end{theorem}

Observe that, with $r=\omega-1$, Theorem \ref{le3mu} gives the following
inequality from \cite{Nik02}; it has been applied to obtain a two line proof
of the spectral precise Tur\'{a}n theorem in \cite{Nik07}.

\begin{theorem}
If $G$ is a graph with $\omega\left(  G\right)  =\omega\geq2$ and $\mu\left(
G\right)  =\mu,$ then
\[
\mu^{\omega}\leq k_{2}\left(  G\right)  \mu^{\omega-2}+2k_{3}\left(  G\right)
\mu^{\omega-3}+\cdots+\left(  \omega-1\right)  k_{\omega}\left(  G\right)  .
\]

\end{theorem}

Another important consequence of Theorem \ref{le3mu}, also in \cite{BoNi07},
gives a lower bound on the number of cliques of any order as stated below.

\begin{theorem}
\label{tmomo}If $r\geq2$ and $G\in\mathcal{G}\left(  n\right)  ,$ then
\[
k_{r+1}\left(  G\right)  \geq\left(  \frac{\mu\left(  G\right)  }{n}%
-1+\frac{1}{r}\right)  \frac{r\left(  r-1\right)  }{r+1}\left(  \frac{n}%
{r}\right)  ^{r+1}.
\]

\end{theorem}

The remaining two theorems of this subsection are given in \cite{Nik07b} and
have multiple uses. The first one relates the numbers of triangles, edges and
vertices of a graph with the smallest eigenvalue of its adjacency matrix.

\begin{theorem}
\label{tmutm}If $G\in\mathcal{G}\left(  n,m\right)  ,$ then%
\begin{equation}
\mu_{n}\left(  G\right)  \leq\frac{3n^{3}k_{3}\left(  G\right)  -4m^{3}%
}{nm\left(  n^{2}-2m\right)  } \label{ineq}%
\end{equation}
with equality if and only if $G$ is a regular complete multipartite graph.
\end{theorem}

Inequality (\ref{ineq}) should be regarded as a multifaceted relation that can
be used for different purposes. By way of illustration, let us restate it as a
lower bound on $k_{3}\left(  G\right)  ,$ getting
\begin{equation}
k_{3}\left(  G\right)  \geq\frac{\mu_{n}\left(  G\right)  \left(  nm\left(
n^{2}-2m\right)  \right)  +4m^{3}}{3n^{3}}, \label{ineq1}%
\end{equation}
with equality holding for regular complete multipartite graphs. However, for
all dense quasi-random graphs we have $\mu_{n}\left(  G\right)  =o\left(
n\right)  $ and $3k_{3}\left(  G\right)  =4\left(  1+o\left(  1\right)
\right)  m^{3}/n^{2}.$ This implies that
\[
\frac{4m^{3}}{3n^{3}}+o\left(  1\right)  \frac{m^{3}}{n^{3}}=k_{3}\left(
G\right)  \geq o\left(  1\right)  mn+\frac{4m^{3}}{3n^{3}},
\]
and we reach the somewhat paradoxical conclusion that inequality (\ref{ineq1})
is tight up to low order additive terms for almost all graphs, since almost
all graphs are dense and quasi-random.

Statements similar to Theorem \ref{tmutm} have been obtained in \cite{Nik07b}
for the largest Laplacian eigenvalue $\lambda_{n}\left(  G\right)  $ as well.

\begin{theorem}
If $G\in\mathcal{G}\left(  n,m\right)  ,$ then%
\[
6nk_{3}\left(  G\right)  \geq\left(  n+\lambda_{n}\left(  G\right)  \right)
\sum_{u\in V\left(  G\right)  }d^{2}\left(  u\right)  -2nm\lambda_{n}\left(
G\right)
\]
with equality if and only if $G$ is a complete multipartite graph, and
\[
\lambda_{n}\left(  G\right)  \geq\frac{2m^{2}-3nk_{3}\left(  G\right)
}{m\left(  n^{2}-2m\right)  }n,
\]
with equality if and only if $G$ is a regular complete multipartite graph.
\end{theorem}

\subsection{Chromatic number}

Let $G$ be a graph of order $n.$ One of the best known results in spectral
graph theory is the inequality of A.J. Hoffman \cite{Hof70}
\begin{equation}
\chi\left(  G\right)  \geq1+\frac{\mu_{1}\left(  G\right)  }{-\mu_{n}\left(
G\right)  }, \label{hofin}%
\end{equation}
However, it seems that there is a lot more to find in this area. Indeed, in
\cite{Nik07a} we proved the following alternative bound.

\begin{theorem}
For every graphs of order $n,$%
\begin{equation}
\chi\left(  G\right)  \geq1+\frac{\mu_{1}\left(  G\right)  }{\lambda
_{n}\left(  G\right)  -\mu_{1}\left(  G\right)  }. \label{mainin}%
\end{equation}
Equality holds if and only if every two color classes of $G$ induce a regular
bipartite graph of degree $\left\vert \mu_{n}\left(  G\right)  \right\vert $.
\end{theorem}

When $G$ is obtained from $K_{n}$ by deleting an edge, inequality
(\ref{mainin}) gives $\chi\left(  G\right)  =n-1,$ while (\ref{hofin}) gives
only $\chi\left(  G\right)  \geq n/2+2.$ By contrast, for a sufficiently large
wheel $W_{1,n}$, i.e., a vertex joined to all vertices of a cycle of length
$n$, (\ref{mainin}) gives $\chi\geq2,$ while (\ref{hofin}) gives $\chi\geq3.$

However, such comparisons are not too informative since, in \cite{Nik07a},
both (\ref{mainin}) and (\ref{hofin}) have been deduced from the same matrix theorem.

\section{\label{tool}Some useful tools}

In this section we present some results that we have found useful on multiple
occasions. The selection and the arrangement of these results does not follow
any particular pattern.

We start with an inequality stated by Moon and Moser in \cite{MoMo62}; it
seems that Khad\v{z}iivanov and Nikiforov \cite{KhNi78} were the first to
publish its complete proof, see also \cite{Lov79}, Problem 11.8. The
inequality has been used in many questions, say in the proof of Theorem
\ref{ESnew}.

\begin{lemma}
\label{momo}Let $1\leq s<t<n$, and let $G$ be a graph of order $n,$ with
$k_{t}\left(  G\right)  >0$. Then
\begin{equation}
\frac{\left(  t+1\right)  k_{t+1}\left(  G\right)  }{tk_{t}\left(  G\right)
}-\frac{n}{t}\geq\frac{\left(  s+1\right)  k_{s+1}\left(  G\right)  }%
{sk_{s}\left(  G\right)  }-\frac{n}{s}. \label{MoMo}%
\end{equation}

\end{lemma}

The following two simple lemmas were used to obtain a number of results in
Section \ref{class}. The first one was proved in \cite{Nik08}, and the second
one in \cite{Nik10}.

\begin{lemma}
\label{routl}Let $r\geq2,$ let $c,n,m,s$ be such that
\[
0<c\leq1/2,\text{ \ \ }n\geq\exp\left(  c^{-r}\right)  ,\text{ \ \ }%
s=\left\lfloor c^{r}\log n\right\rfloor \leq\left(  c/2\right)  m+1,
\]
and let $G$ be a bipartite graph with parts $A$ and $B$ of size $m$ and $n.$
If $e\left(  G\right)  \geq cmn,$ then $G$ contains a $K_{2}\left(
s,t\right)  $ with parts $S\subset A$ and $T\subset B$ such that $\left\vert
S\right\vert =s$ and $\left\vert T\right\vert =t>n^{1-c^{r-1}}$.
\end{lemma}

\begin{lemma}
Let $\alpha,c,n,m$ be such that
\[
0<\alpha\leq1,\text{ \ \ }1\leq c\log n\leq\alpha m/2+1,
\]
and let $G$ be a bipartite graph with parts $A$ and $B$ of size $m$ and $n.$
If $e\left(  G\right)  \geq\alpha mn,$ then $G$ contains a $K_{2}\left(
s,t\right)  $ with parts $S\subset A$ and $T\subset B$ such that $\left\vert
S\right\vert =\left\lfloor c\log n\right\rfloor $ and $\left\vert T\right\vert
=t>n^{1-c\log\alpha/2}$.
\end{lemma}

The following lemma, given in \cite{Nik09f}, strengthens a classical condition
for the existence of paths given by Erd\H{o}s and Gallai \cite{ErGa59}. It has
been used to obtain results about forbidden cycles and elsewhere.

\begin{lemma}
\label{lev6}Suppose that $k\geq1$ and let the vertices of a graph $G$ be
partitioned into two sets $U$ and $W$.

(i) If
\[
2e\left(  U\right)  +e\left(  U,W\right)  >\left(  2k-2\right)  \left\vert
U\right\vert +k\left\vert W\right\vert ,
\]
then there exists a path of order $2k$ or $2k+1$ with both ends in $U.$

(ii) If
\[
2e\left(  U\right)  +e\left(  U,W\right)  >\left(  2k-1\right)  \left\vert
U\right\vert +k\left\vert W\right\vert ,
\]
then there exists a path of order $2k+1$ with both ends in $U.$
\end{lemma}

The following lemma from \cite{Nik10c} was used to prove Theorems \ref{extBT},
\ref{extJin} and \ref{extCJK}, but may be used to carry over other stability
results from triangle-free graphs to $K_{r}$-free graphs for $r>3.$

\begin{lemma}
\label{reduL}Let $r\geq3$ and let $G$ be a maximal $K_{r+1}$-free graph of
order $n.$ If
\[
\delta\left(  G\right)  >\left(  1-\frac{2}{2r-1}\right)  n,
\]
then $G$ has a vertex $u$ such that the vertices not joined to $u$ are independent.
\end{lemma}

The following lemma, given in \cite{Nik09g}, bounds the minimal entry of
eigenvectors to the spectral radius of the adjacency matrix. This can be
useful in various situations, e.g., in conjunction with Lemma \ref{lev2} from
\cite{Nik10a} and Theorem \ref{thv4} it can be used to prove upper bounds on
$\mu\left(  G\right)  $ by induction. Both lemmas have been used to prove
several results in Section \ref{spec}.

\begin{lemma}
\label{lev1}Let $G$ be a graph of order $n$ with minimum degree $\delta\left(
G\right)  =\delta$ and $\mu\left(  G\right)  =\mu.$ If $\left(  x_{1}%
,\ldots,x_{n}\right)  $ is a unit eigenvector to $\mu,$ then%
\[
\min\left\{  x_{1},\ldots,x_{n}\right\}  \leq\sqrt{\frac{\delta}{\mu
^{2}+\delta n-\delta^{2}}}.
\]

\end{lemma}

\begin{lemma}
\label{lev2}Let $G$ be a graph of order $n$ and let $\left(  x_{1}%
,\ldots,x_{n}\right)  $ be a unit eigenvector to $\mu\left(  G\right)  .$ If
$u$ is\ a vertex satisfying $x_{u}=\min\left\{  x_{1},\ldots,x_{n}\right\}  ,$
then%
\[
\mu\left(  G-u\right)  \geq\mu\left(  G\right)  \frac{1-2x_{u}^{2}}%
{1-x_{u}^{2}}.
\]

\end{lemma}

The theorem below, given in \cite{Nik09g}, has been used to prove the spectral
analog of several nonspectral results.

\begin{theorem}
\label{thv4}Let $\alpha,\beta,\gamma,K$ and $n$ be such that
\[
0<4\alpha\leq1,\text{ \ \ }0<2\beta\leq1,\text{ \ \ }1/2-\alpha/4\leq
\gamma<1,\text{ \ \ }K\geq0,\text{ \ \ }n\geq\left(  42K+4\right)  /\alpha
^{2}\beta,
\]
and let $G$ be a graph of order $n.$ If
\[
\mu\left(  G\right)  >\gamma n-K/n\text{ \ \ and \ \ }\delta\left(  G\right)
\leq\left(  \gamma-\alpha\right)  n,
\]
then there exists an induced subgraph $H\subset G\ $with $\left\vert
H\right\vert \geq\left(  1-\beta\right)  n,$ satisfying one of the following conditions:

(a) $\mu\left(  H\right)  >\gamma\left(  1+\beta\alpha/2\right)  \left\vert
H\right\vert ;$

(b) $\mu\left(  H\right)  >\gamma\left\vert H\right\vert $ and $\delta\left(
H\right)  >\left(  \gamma-\alpha\right)  \left\vert H\right\vert .$
\end{theorem}

The abundance of parameters in Theorem \ref{thv4} may obstruct its
understanding. In summary, the theorem can be applied when one has to prove
that if $\mu\left(  G\right)  $ is sufficiently large then $G$ contains some
subgraph $F.$ If $\delta\left(  G\right)  $ is not large enough, by tossing
away not too many low degree vertices, one gets a graph $H$ in which either
both $\mu\left(  H\right)  $ and $\delta\left(  H\right)  $ are large enough
or $\mu\left(  H\right)  $ is considerably above the expected average. Most
likely, either of these properties will help to find a copy of $F$ in $H.$ The
many parameters ensure greater flexibility.

In \cite{BoNi04a}, using interlacing, Bollob\'{a}s and Nikiforov gave the
following inequality, which has been used to prove several results involving
the minimum eigenvalue of the adjacency matrix, e.g., Theorem \ref{minK}.

\begin{theorem}
\label{BoNi}If $G\in\mathcal{G}\left(  n,m\right)  ,$ then for every partition
$V\left(  G\right)  =V_{1}\cup V_{2},$%
\[
\mu_{n}\left(  G\right)  \leq\frac{2e\left(  V_{1}\right)  }{\left\vert
V_{1}\right\vert }+\frac{2e\left(  V_{2}\right)  }{\left\vert V_{2}\right\vert
}-\frac{2m}{n}.
\]

\end{theorem}

Note that this inequality is analogous to the well-known inequality for the
Laplacian (see Mohar, \cite{Moh96}):%
\[
\lambda_{n}\left(  G\right)  \geq\frac{e\left(  V_{1},V_{2}\right)
}{\left\vert V_{1}\right\vert \left\vert V_{2}\right\vert }n,
\]
and in fact for regular graphs both inequalities are identical.

\section{\label{spf}Illustration proofs}

The purpose of this section is to illustrate the use of the tools developed
for translating nonspectral into spectral results. To this end we shall sketch
the proofs of Theorems \ref{tStab} and \ref{stStab}.

The structure of both proofs is identical. In both proofs we shall use Theorem
\ref{Nik1} from Section \ref{ErSt}. The main difference comes from the fact
that in the proof of Theorem \ref{tStab} we use Theorem \ref{stabj} while in
the proof of Theorems \ref{stStab} we use the analogous spectral result
Theorem \ref{spN1.2}.\bigskip

\textbf{Proof of Theorem \ref{tStab} }Let $G$ be a graph of order $n$ with
$e\left(  G\right)  >\left(  1-1/r-\varepsilon\right)  n^{2}/2.$ Define the
procedure $\mathcal{P}$ as follows:\medskip

\textbf{While}\emph{ }$js_{r+1}\left(  G\right)  >n^{r-1}/r^{r+6}$ \textbf{do}

\qquad\emph{Select an edge contained in }$\left\lceil n^{r-1}/r^{r+6}%
\right\rceil $\emph{ cliques of order }$r+1$ \emph{and remove it from }%
$G.$\medskip

Set for short $\theta=c^{1/\left(  r+1\right)  }r^{r+6}$ and assume first that
$\mathcal{P}$ removes at least $\left\lceil \theta n^{2}\right\rceil $ edges
before stopping. Then
\[
k_{r+1}\left(  G\right)  \geq\theta n^{r-1}/r^{r+6}=c^{1/\left(  r+1\right)
}n^{r+1},
\]
and Theorem \ref{Nik1} implies that
\[
K_{r+1}\left(  \left\lfloor c\ln n\right\rfloor ,\ldots,\left\lfloor c\ln
n\right\rfloor ,\left\lceil n^{1-\sqrt{c}}\right\rceil \right)  \subset G.
\]
Thus, in this case condition \emph{(a)} holds, completing the proof.

Assume therefore that $\mathcal{P}$ removes fewer than $\left\lceil \theta
n^{2}\right\rceil $ edges before stopping. Writing $G^{\prime}$ for the
resulting graph, we see that
\[
e\left(  G^{\prime}\right)  >e\left(  G\right)  -\theta n^{2}>\left(
1-1/r-\varepsilon-\theta\right)  n^{2}/2
\]
and $js_{r+1}\left(  G^{\prime}\right)  <n^{r-1}/r^{r+6}.$ Here we want to
apply Theorem \ref{stabj} and so we check for its prerequisites. First, from
$\log n\geq1/c\geq r^{3\left(  r+14\right)  \left(  r+1\right)  }$ we easily
get $n>r^{8}.$ Also,
\[
\varepsilon+\theta<r^{-8}/8.
\]
Now, Theorem \ref{stabj} implies that $G^{\prime}$ contains an induced
$r$-partite subgraph $G_{0}$ satisfying
\[
\left\vert G_{0}\right\vert \geq\left(  1-\sqrt{2\left(  \varepsilon
+\theta\right)  }\right)  n\text{ and }\delta\left(  G_{0}\right)  >\left(
1-1/r-2\sqrt{2\left(  \varepsilon+\theta\right)  }\right)  n.
\]

By routine calculations we find that $G$ differs from $T_{r}\left(  n\right)
$ in fewer than
\[
\left(  \theta+\left(  2r^{2}-r\right)  \sqrt{2\left(  \varepsilon
+\theta\right)  }\right)  n^{2}<\left(  \varepsilon^{1/3}+c^{1/\left(
3r+3\right)  }\right)  n^{2}%
\]
edges, and condition \emph{(b)} follows, completing the proof of Theorem
\ref{tStab}.$\hfill\square$

\bigskip

\textbf{Proof of Theorem \ref{stStab} }Let $G$ be a graph of order $n$ with
$\mu\left(  G\right)  >\left(  1-1/r-\varepsilon\right)  n.$ Define the
procedure $\mathcal{P}$ as follows:\medskip

\textbf{While}\emph{ }$js_{r+1}\left(  G\right)  >n^{r-1}/r^{2r+5}$
\textbf{do}

\qquad\emph{Select an edge contained in }$\left\lceil n^{r-1}/r^{2r+5}%
\right\rceil $\emph{ cliques of order }$r+1$ \emph{and remove it from }%
$G.$\medskip

Set for short $\theta=c^{1/\left(  r+1\right)  }r^{2r+5}$ and assume first
that $\mathcal{P}$ removes at least $\left\lceil \theta n^{2}\right\rceil $
edges before stopping. Then
\[
k_{r+1}\left(  G\right)  \geq\theta n^{r-1}/r^{2r+5}=c^{1/\left(  r+1\right)
}n^{r+1},
\]
and Theorem \ref{Nik1} implies that
\[
K_{r+1}\left(  \left\lfloor c\ln n\right\rfloor ,\ldots,\left\lfloor c\ln
n\right\rfloor ,\left\lceil n^{1-\sqrt{c}}\right\rceil \right)  \subset G.
\]
Thus, in this case condition \emph{(a)} holds, completing the proof.

Assume now that $\mathcal{P}$ removes fewer than $\left\lceil \theta
n^{2}\right\rceil $ edges before stopping. Write $G^{\prime}$ for the
resulting graph; we obviously have $js_{r+1}\left(  G^{\prime}\right)  \leq
n^{r-1}/r^{2r+5}.$ Letting $\mu\left(  X\right)  $ be the largest eigenvalue
of a Hermitian matrix $X,$ recall Weyl's inequality
\[
\mu\left(  B\right)  \geq\mu\left(  A\right)  -\mu\left(  A-B\right)  ,
\]
holding for any Hermitian matrices $A$ and $B.$ Also, recall that $\mu\left(
H\right)  \leq\sqrt{2e\left(  H\right)  }$ for any graph $H.$ Applying these
results to the graphs $G$ and $G^{\prime},$ we find that
\[
\mu\left(  G^{\prime}\right)  \geq\mu\left(  G\right)  -\sqrt{2\theta}%
n\geq\left(  1-1/r-\varepsilon-\sqrt{2\theta}\right)  n.
\]
Here we want to apply Theorem \ref{spN1.2} and so we check for its
prerequisites. First, from $\log n\geq1/c\geq r^{8\left(  r+21\right)  \left(
r+1\right)  }$ we easily get $n>r^{20}.$ Also,
\[
\varepsilon+\sqrt{2\theta}<2^{-10}r^{-6}.
\]
Now, Theorem \ref{spN1.2} implies that $G^{\prime}$ contains an induced
$r$-partite subgraph $G_{0},$ satisfying
\[
\left\vert G_{0}\right\vert \geq\left(  1-4\left(  \varepsilon+\sqrt{2\theta
}\right)  ^{1/3}\right)  n\text{ and }\delta\left(  G_{0}\right)  >\left(
1-1/r-7\left(  \varepsilon+\sqrt{2\theta}\right)  ^{1/3}\right)  n.
\]

By routine calculations we find that $G$ differs from $T_{r}\left(  n\right)
$ in fewer than
\[
\left(  \theta+\left(  7r^{2}-3r\right)  \left(  \varepsilon+\sqrt{2\theta
}\right)  ^{1/3}\right)  n^{2}<\left(  \varepsilon^{1/4}+c^{1/\left(
8r+8\right)  }\right)  n^{2}%
\]
edges, and condition \emph{(b)} follows, completing the proof of Theorem
\ref{stStab}.$\hfill\square$

\section{\label{not}Notation and basic facts}

Throughout the survey our notation generally follows \cite{Bol98}. Given a
graph $G,$ we write:\medskip

- $V\left(  G\right)  $ for the vertex set of $G;$

- $E\left(  G\right)  $ for the edge set of $G$ and $e\left(  G\right)  $ for
$\left\vert E\left(  G\right)  \right\vert ;$

- $\alpha\left(  G\right)  $ for the independence number of $G$ (see below);

- $\delta\left(  G\right)  $ and $\Delta\left(  G\right)  $ for the minimum
and maximum degrees of $G;$

- $\omega\left(  G\right)  $ for the clique number of $G$ (see below);

- $k_{s}\left(  G\right)  $ for the number of $s$-cliques of $G$ (see below);

- $G-u$ for the graph obtained by removing the vertex $u\in V\left(  G\right)
;$

- $\Gamma\left(  u\right)  $ for the set of neighbors of a vertex $u,$ and
$d\left(  u\right)  $ for $\left\vert \Gamma\left(  u\right)  \right\vert ;$

- $e\left(  X\right)  $ for the number of edges induced by a set $X\subset
V\left(  G\right)  ;$

- $e\left(  X,Y\right)  $ for the number of edges joining vertices in $X$ to
vertices in $Y,$ where $X$ and $Y$ are disjoint subsets of $V\left(  G\right)
;$\medskip

We write $\mathcal{G}(n)$ for the set of graphs of order $n$ and
$\mathcal{G}\left(  n,m\right)  $ for the set of graphs of order $n$ and size
$m.$

Also, $\left[  n\right]  $ stands for the set $\left\{  1,2,\ldots,n\right\}
.$\bigskip

\textbf{Mini glossary}\medskip

\textbf{clique -} a subgraph that is complete. An $s$-clique has $s$ vertices;
$k_{s}\left(  G\right)  $ stands for the number of $s$-cliques of $G;$

\textbf{clique number -} the size of the largest clique of $G,$ denoted by
$\omega\left(  G\right)  ;$

\textbf{chromatic number }- the minimum number of independent sets that
partition $V\left(  G\right)  ,$ denoted by $\chi\left(  G\right)  ;$

\textbf{disjoint union }of two graphs $G$ and $H$ is the union of two vertex
disjoint copies of $G$ and $H.$ The disjoint union of $G$ and $H$ is denoted
by $G+H;$

\textbf{independent set -} a set of vertices of $G$ that induces no edges;

\textbf{independence number -} the size of the largest independent set of $G,$
denoted by $\alpha\left(  G\right)  ;$

\textbf{join }of two vertex disjoint graphs $G$ and $H$ is the union of $G$
and $H$ together with all edges between $G$ and $H.$ The join of $G$ and $H$
is denoted by $G\vee H;$

\textbf{joint }- a set of cliques of the same order sharing an edge. An
$r$-joint of size $t$ consists of $t$ cliques of order $r;$

\textbf{book of size }$t$ - a $3$-joint of size $t,$ that is to say, a
collection of $t$ triangles sharing an edge;

\textbf{homomorphic graph - }a graph $G$ is said to be homomorphic to a graph
$H,$ if there exists a map $f:V\left(  G\right)  \rightarrow V\left(
H\right)  $ such that $uv\in E\left(  G\right)  $ implies $f\left(  u\right)
f\left(  v\right)  \in E\left(  H\right)  ;$

\textbf{graph property -} a family of graphs closed under isomorphisms;

\textbf{hereditary property - }graph property closed under taking induced subgraphs;

\textbf{monotone property - }graph property closed under taking subgraphs;

$H$\textbf{-free graph: }a graph that has no subgraph isomorphic to $H;$

\textbf{friendship graph - }a collection of triangles sharing a single common vertex;

$k$-\textbf{th power of a cycle} $C_{n}$ - a graph with vertices $\left\{
1,2,\ldots,n\right\}  ,$ and $\left(  i,j\right)  $ is an edge if
$i-j=\pm1,\pm2,\cdots,\pm k$ $\operatorname{mod}$ $n$;

$K_{r}$ and $\overline{K}_{r}$ - the complete and the edgeless graph of order
$r;$

$K_{r}\left(  s_{1},s_{2},...,s_{r}\right)  $ - the complete $r$-partite graph
with class sizes $s_{1},s_{2},...,s_{r}.$ We set for short%
\[
K_{r}\left(  p\right)  =K_{r}\left(  p,...,p\right)  \text{ \ \ and \ \ }%
K_{r}\left(  p;q\right)  =K_{r}\left(  p,...,p,q\right)  ;
\]

$r$-\textbf{uniform hypergraph - }a hypergraph whose edges are subsets of $r$ vertices;

\textbf{Tur\'{a}n graph }$T_{r}\left(  n\right)  $ - given $n\geq r\geq2,$
this is the complete $r$-partite graph whose class sizes differ by at most
one. We let $t_{r}\left(  n\right)  =e\left(  T_{r}\left(  n\right)  \right)
$. If $t$ is the remainder of $n$ $\operatorname{mod}$ $r,$ then
\[
t_{r}\left(  n\right)  =\frac{r-1}{2r}\left(  n^{2}-t^{2}\right)  +\binom
{t}{2},
\]
which in turn implies that%
\[
\frac{r-1}{2r}n^{2}-\frac{r}{8}\leq t_{r}\left(  n\right)  \leq\frac{r-1}%
{2r}n^{2};
\]

\textbf{Tur\'{a}n problem }- given a family of graphs $F,$ find the maximum
number of edges in a graph of order $n$, having no subgraph belonging to $F$;

\textbf{quasi-random graph} - informally, an almost regular graph, in which
the second largest in modulus eigenvalue is much smaller than the spectral radius;

\textbf{spectral radius of a graph - }in general, the spectral radius of a
matrix is the largest modulus of its eigenvalues. For a graph, this is usually
the spectral radius of its adjacency matrix, which is an eigenvalue itself;

\textbf{Laplacian matrix - }the matrix $L=D-A,$ where $A$ is the adjacency
matrix and $D$ is the diagonal matrix of the row-sums of $A,$ that is the
degrees of $G;$

$Q$\textbf{-matrix, }also known as \textbf{signless Laplacian - }the matrix
$Q=D+A;$

\textbf{Szemer\'{e}di's Regularity Lemma} - an important result of analytical
graph theory, which states that every graph can be approximated by graphs of
bounded order. For background on this lemma we refer the reader to
\cite{Bol98}, Section IV.5;

\textbf{Zarankiewicz problem - }a class of problems aiming to determine the
maximum number of edges in a graph with no $K_{s,t}.$ There are several
variations, most of which are only partially solved. See \cite{Bol98} for
details.\bigskip

\textbf{Acknowledgement }I am most grateful to the referee for the efficient
and kind help.

\end{document}